\documentclass{amsart}

\usepackage{breqn}
\usepackage{mathtools}
\usepackage{accents}

\begin{document}

\title[curvature]{Curvature on the integers, II}
\bigskip

\def \h{\hat{\ }}
\def \cO{\mathcal O}
\def \ra{\rightarrow}
\def \bZ{{\mathbb Z}}
\def \cP{\mathcal V}
\def \cH{{\mathcal H}}
\def \cB{{\mathcal B}}
\def \d{\delta}
\def \cC{{\mathcal C}}
\def \jor{\text{jor}}

\newtheorem{THM}{{\!}}[section]
\newtheorem{THMX}{{\!}}
\renewcommand{\theTHMX}{}
\newtheorem{theorem}{Theorem}[section]
\newtheorem{corollary}[theorem]{Corollary}
\newtheorem{lemma}[theorem]{Lemma}
\newtheorem{proposition}[theorem]{Proposition}
\theoremstyle{definition}
\newtheorem{definition}[theorem]{Definition}
\theoremstyle{remark}
\newtheorem{remark}[theorem]{Remark}
\newtheorem{example}[theorem]{\bf Example}
\numberwithin{equation}{section}
\address{Department of Mathematics and Statistics\\University of New Mexico \\ Albuquerque, NM 87131, USA}
\email{buium@math.unm.edu} 
\maketitle

\bigskip

\medskip
\centerline{\bf Alexandru Buium}
\bigskip

\begin{abstract} 
In a prequel to this paper \cite{curvature1} a notion of curvature on the integers was introduced, based on the technique of  ``analytic continuation between primes", introduced in \cite{laplace}. In this paper, which is essentially independent of its prequel, we introduce another notion of curvature on the integers, based on ``algebraization of Frobenius lifts by correspondences." Our main results are vanishing/non-vanishing theorems for this new type of curvature in the case of ``Chern connections" attached to classical groups.  
\end{abstract}

\section{Introduction}

This paper is, in principle, a continuation of \cite{curvature1} but, from a logical standpoint, it is independent of \cite{curvature1}.  For the motivation of our theory, and its comparison with classical differential geometry, we refer to the discussion in \cite{curvature1}. More generally, the present paper should be viewed as taking a step in the direction of developing a ``differential geometry on $Spec\ \bZ$"; this direction of research is consistent with the study in  \cite{char, book, adel2, curvature1} of ``arithmetic differential equations," as well as with Borger's viewpoint in \cite{borgerf1} on the ``field with one element."
 
In \cite{curvature1} we started by viewing the ring of integers,
  $\bZ$, as an analogue of a ring of functions on an infinite dimensional manifold in which the various directions are given by the primes; then, in the spirit of \cite{char, book, laplace}, we replaced the  partial derivative operators, acting on functions on a manifold,  by  Fermat quotient type  operators, called {\it $p$-derivations}, acting on numbers.  We then developed an arithmetic analogue  of connections and curvature on the ``manifold" $Spec\ \bZ$ and we proved a series of vanishing/non-vanishing results for the curvature of ``Chern connections" attached to the classical groups.   In order to achieve this program we had to deal, in \cite{curvature1}, with the following  difficulty: the various $p$-derivations defining the Chern  connections on $GL_n$  are defined as self-maps  of the corresponding $p$-adic completions of the ring of functions of $GL_n$ so, when $p$ varies, the $p$-derivations under consideration  do not act on the same ring. Consequently, one cannot directly consider their commutator and, hence, their curvature.
  In \cite{curvature1} we overcame this problem by implementing the technique of {\it analytic continuation between primes} introduced in \cite{laplace}; this technique only worked in the case of classical groups defined by symmetric/antisymmetric matrices with entries roots of unity or zero.
  In the present paper we will overcome the above mentioned difficulty in a different way, namely by ``algebraizing" the analytic picture in \cite{curvature1}. 
  This algebraization method has at least two  advantages: 1)  it works for arbitrary symmetric/antisymmetric matrices, with entries not necessarily roots of unity or zero and 2) it deals with schemes, indeed with function fields of varieties, rather than with  formal schemes.  The price to pay is that one needs to replace endomorphisms by correspondences. On the other hand correspondences can be composed and commutators can be attached to them, leading to a new concept of curvature.
   The resulting picture, in the present paper, will then acquire, as we shall see,  a  ``birational/motivic" flavor. Our main results will, again, be vanishing/non-vanishing theorems for (this new type of) curvature in the case of  ``Chern connections" attached to the classical groups.  
   Our {\it Chern connections} are analogues of the Chern connections on hermitan vector bundles \cite{GH, kobayashi} and were introduced in \cite{adel2}; we will review their definition following \cite{adel2} presently.

In the rest of this Introduction we give a rough idea of our main constructions and results.  We begin by recalling from the Introduction to \cite{curvature1} a few basic definitions. 
Recall from \cite{char, joyal} that  a $p$-{\it derivation} on a $p$-torsion free ring $B$  is a set theoretic map $\d_p:B\ra B$ such that the map
$\phi_p:B\ra B$ defined by $\phi_p(b):=b^p+p\d_p (b)$ is a ring homomorphism; we say $\phi_p$ is the {\it lift of Frobenius} attached to $\d_p$.
Throughout the paper  $A$ will denote  the ring $\bZ[1/M,\zeta_N]$ where $M$ is some even integer and $\zeta_N$ is a primitive $N$-th root of unity, $N\geq 1$.  Also we let $G=GL_n=Spec\ A[x,\det(x)^{-1}]$ be the general linear group scheme over $A$, where $x=(x_{kl})$ is an $n \times n$ matrix of ideterminates. 
Fix a (possibly infinite) set ${\mathcal V}$ of prime integers $p$ not dividing $MN$.
By an {\it adelic connection} on $G$ we understand  a family $(\d_p)$, indexed by  ${\mathcal V}$, where, 
for each $p$, $\d_p$ is 
a $p$-derivation on the $p$-adic completion $A[x,\det(x)^{-1}]^{\widehat{p}}$ of $A[x, \det(x)^{-1}]$. We can consider the  attached  family $(\phi_p)$ of lifts of Frobenius 
on the rings 
$A[x,\det(x)^{-1}]^{\widehat{p}}$;
 we shall identify the $\phi_p$'s with endomorphisms of the 
   $p$-adic completion $G^{\widehat{p}}$ of $G$. Recall that there is a bijection between the set of adelic connections $(\d_p)$ on $G$ and the set of families $(\Delta_p)$ where, for each $p$,  $\Delta_p$ is an $n\times n$ matrix with entries in $A[x,\det(x)^{-1}]^{\widehat{p}}$; the bijection is provided by $\d_p x=\Delta_p$ and we have $\phi_p(x)=x^{(p)}+p\Delta_p$, where $x^{(p)}:=(x_{kl}^p)$.

Natural examples of adelic connections were introduced in  \cite{adel2} as follows.
Let $q\in GL_n(A)$ with $q^t=\pm q$, where the $t$ superscript means ``transpose". (Morally if $q^t=q$ then $q$ can be viewed as an analogue of a metric on a principal bundle over $Spec\ \bZ$; the case $q^t=-q$ corresponds to $2$-forms; cf. \cite{curvature1}.)
 One can attach  to $q$  maps $\cH_q:G\ra G$ and  $\cB_q:G\times G\ra G$ defined by
$\cH_q(x)=x^tqx$, $\cB_q(x,y)=x^tqy$.
We continue to denote by $\cH_q,\cB_q$ the maps induced on the $p$-adic completions $G^{\widehat{p}}$ and $G^{\widehat{p}}\times G^{\widehat{p}}$. One can consider, in addition to the data above, the adelic connection $\d_0=(\d_{0p})$ on $G$ with $\d_{0p}x=0$. (Morally $\d_0$ is thought of as fixing a $\overline{\partial}$ operator, or a complex structure; cf. \cite{curvature1}.) Denote by $(\phi_p)$ and $(\phi_{0p})$ the families of lifts of Frobenius attached to $\d$ and $\d_0$ respectively.
In \cite{adel2} it was shown that for any $q$ as above there exists a unique adelic connection $\d$ such that the following diagrams are commutative:
\begin{equation}
 \label{coconut}
  \begin{array}{rcl}
G^{\widehat{p}} & \stackrel{\phi_p}{\longrightarrow} & G^{\widehat{p}}\\
\cH_q  \downarrow & \  & \downarrow \cH_q \\
G^{\widehat{p}} & \stackrel{\phi_{0p}}{\longrightarrow} & G^{\widehat{p}}\\
\end{array}\ \ \ \ \ \ 
\begin{array}{rcl}
G^{\widehat{p}} & \stackrel{\phi_{0p} \times \phi_p}{\longrightarrow} & G^{\widehat{p}}\times G^{\widehat{p}}\\
\phi_p \times \phi_{0p} \downarrow & \  & \downarrow \cB_q\\
 G^{\widehat{p}}\times G^{\widehat{p}} & \stackrel{\cB_q}{\longrightarrow} & G^{\widehat{p}}\end{array}\end{equation}
 In \cite{curvature1} we  called $\d=(\d_p)$ the {\it Chern connection} attached to $q$; for the analogy with classical differential geometry \cite{GH, kobayashi} see \cite{curvature1}.
 
  In what follows we would like to introduce curvature of adelic connections via a construction involving correspondences. So we need to introduce some terminology related to correspondences. 
 Let $E$ be the fraction field of $G=GL_n$; so, if $K={\mathbb Q}(\zeta_N)$ is the fraction field of $A$, then  $E=K(x)$ is a purely transcendental extension of $K$ generated by the variables $x_{kl}$.
 By a {\it correspondence} on $E$ we will understand a triple $\Gamma=(Y,\pi,\varphi)$ where $Y$ is a reduced  non-empty scheme and $\pi,\varphi:Y\ra Spec\ E$ are finite  morphisms
  of schemes. So $Y=Spec\ F$ is the spectrum of a finite product, $F$, of fields each of which is a finite extension of $E$ via both $\pi$ and $\varphi$; the degrees of these two maps are referred to as the {\it left} and {\it right} degree of $\Gamma$ respectively. If $Y$ is irreducible, i.e., the spectrum of a field, we say $\Gamma$ is {\it irreducible}.
  For any such correspondence $\Gamma=(Y,\pi,\varphi)$ we still denote by 
  $\pi,\varphi:E\ra F$ the induced ring homomorphisms, we let $\text{tr}_{\pi}:F\ra E$ be the trace map corresponding to the morphism $\pi$, and
 we denote by $\Gamma^*:E\ra E$  the additive group homomorphism $\Gamma^*=\text{tr}_{\pi}\circ \varphi:E\ra F\ra E$. 
 By a {\it correspondence structure} for  an adelic connection $\d=(\d_p)$ on  $G$ we will understand a collection $(\Gamma_p)$ of correspondences on $E$ 
  such that, for each $p$, $\Gamma_p=(Y_p,\pi_p,\varphi_p)$ is {\it  compatible}  with the corresponding lift of Frobenius $\phi_p$ in the following sense: there are morphisms of affine schemes $\pi_{p/G}, \varphi_{p/G}:Y_{p/G}\ra G$, such that  
   $\Gamma_{p}\simeq \Gamma_{p/G}\otimes E$,  $\pi_{p/G}$ is  \'{e}tale, $\pi_{p/G}^{\widehat{p}}:Y_{p/G}^{\widehat{p}}  \ra G^{\widehat{p}}$ is an isomorphism, and $\varphi_{p/G}^{\widehat{p}}=
\phi_p\circ \pi_{p/G}^{\widehat{p}}:Y_{p/G}^{\widehat{p}}\ra G^{\widehat{p}}$.
   Correspondence structures are not a priori unique but all such structures, for a given $\d$, are  {\it compatible} among themselves at all $p$'s that are inert in $K$ (in a sense that will be made precise later; cf. Definition \ref{compatiblel}, Lemma \ref{abovel}, and Remark \ref{abovell}).

 Now, given a correspondence structure $(\Gamma_p)$ on an adelic connection $\d=(\d_p)$ on $GL_n$, we can define the {\it curvature} of $(\Gamma_p)$ as the family $(\varphi_{pp'}^*)$ where $\varphi^*_{pp'}$ is the additive group endomorphism
  \begin{equation}
  \label{varfistar1}
  \varphi^*_{pp'}:=\frac{1}{pp'}(\Gamma^*_{p'}\circ \Gamma^*_p-\Gamma^*_p\circ \Gamma^*_{p'}):E\ra E.
  \end{equation}
  Given one more adelic connection $\overline{\d}=(\overline{\d}_p)=:(\d_{\overline{p}})$ with correspondence structure $(\overline{\Gamma}_p)=:(\Gamma_{\overline{p}})$ one can define the {\it $(1,1)$-curvature} of $(\Gamma_p)$ with respect to $(\Gamma_{\overline{p}})$ as the family $(\varphi_{p\overline{p}'}^*)$ where $\varphi^*_{p\overline{p}'}$ is the additive group endomorphism
  \begin{equation}
  \label{varfistar2}
  \varphi^*_{p\overline{p}'}:=\frac{1}{pp'}(\Gamma^*_{\overline{p}'}\circ \Gamma^*_p-\Gamma^*_p\circ \Gamma^*_{\overline{p}'}):E\ra E\ \ \text{for $p\neq p'$, and}
  \end{equation}
   \begin{equation}
  \label{varfistar22}
  \varphi^*_{p\overline{p}}:=\frac{1}{p}(\Gamma^*_{\overline{p}}\circ \Gamma^*_p-\Gamma^*_p\circ \Gamma^*_{\overline{p}}):E\ra E.
  \end{equation}
  The factors $\frac{1}{pp'}$, $\frac{1}{p}$ are introduced in order to match the definitions in \cite{curvature1} and will play no role in what follows.
  Also the above ``upper $*$" curvatures have a ``lower $*$" version that will be discussed in the body of the paper. In what follows we  let $\overline{\d}$ be equal to $\d_0=(\d_{0p})$, where $\d_{0p}x=0$; we give $\overline{\d}$ the correspondence structure 
$(\Gamma_{\overline{p}})=(Y_{\overline{p}},\pi_{\overline{p}},\varphi_{\overline{p}})$ 
with $Y_{\overline{p}}=Spec\ E$, $\pi_{\overline{p}}$ the identity, and $\varphi_{\overline{p}}(x)=x^{(p)}$. 
We will prove the following results;   more complete results will be proved in the body of the paper. 

 \begin{theorem}
 \label{alta1}
 For any $q\in GL_n(A)$ with $q^t=\pm q$ the Chern connection  on $GL_n$ attached to $q$  admits a correspondence structure.\end{theorem}
 
 Theorem \ref{alta1} takes a very simple (yet non-trivial) form in case $n=1$, $q\in \bZ[1/M]^{\times}\subset A^{\times}=GL_1(A)$; cf. \cite{adel2, curvature1}. Indeed, by loc. cit., in this case, the Chern connection attached to $q$ is given by
 $\phi_p:A[x,x^{-1}]^{\widehat{p}}\ra A[x,x^{-1}]^{\widehat{p}}$,
 \begin{equation}
 \label{formula}
 \phi_p(x)=\left(\frac{q}{p}\right)\cdot q^{\frac{p-1}{2}}\cdot x^p,
 \end{equation}
 where $\left(\frac{q}{p}\right)$ is the Legendre symbol; hence, since $\phi_p(x)\in A[x,x^{-1}]$,  a correspondence structure $(Y_p,\pi_p,\varphi_p)$ of left degree $1$ for this adelic connection can be  given by letting $Y=Spec\ F_p$, $F_p=E$, $\pi_p=id$,
 and $\varphi:E\ra F_p=E$ given by $\varphi_p(x)=\phi_p(x)$ as in \ref{formula}.
 By the way, if $q$ is a square in $A=\bZ[1/M,\zeta_N]$, with one of the square roots $\sqrt{q}\in A^{\times}$ then
 \begin{equation}
  \label{formulala}
 \phi_p(x)=\varphi_p(x)=\frac{(\sqrt{q}\cdot x)^p}{\phi_p(\sqrt{q})}.
 \end{equation}
 By contrast with the above, for $n\geq 2$, $\phi_p(x)$ does not have entries in $A[x\det(x)^{-1}]$ but rather in $A[x\det(x)^{-1}]^{\widehat{p}}$  and correspondences of left degree $>1$ are necessary to algebraize our $\phi_p$'s.
 It would be interesting to  understand what the minimum left degree for the correspondence structures in Theorem \ref{alta1} can be; also one would like to understand if ``algebraic irrational  analogues" of the formula \ref{formulala} can hold for $n\geq 2$.
  
 Note that, in particular,  Theorem \ref{alta1}  allows us to talk about  curvature for the Chern connection attached to any symmetric/antisymmetric $q\in GL_n(A)$, whereas \cite{curvature1} only attaches curvature to Chern connections attached to   $q$'s that have  entries roots of unity or zero. One can ask, however, how the curvature in the present paper compares with the curvature in \cite{curvature1}, in case we are looking at $q$'s with entries roots of unity or zero. So let us
consider the matrices
$$ \left(\begin{array}{rr} 0 & 1\\  - 1 & 0\end{array}\right),\ \ \ 
\left(\begin{array}{rr} 0 & 1\\  1 & 0\end{array}\right)$$
which we refer to as the {\it split}  matrices in $GL_2$ with $q^t=-q$ and $q^t=q$ respectively; they ``define" the ``split" groups $Sp_2$ and $SO_2$ respectively, and were our basic examples in \cite{curvature1} in the $2 \times 2$ case. It turns out that, for these matrices, we  
 have the following results which are parallel to the corresponding results in \cite{curvature1}.

  \begin{theorem}
  \label{alta2}
  Let $q$ be the split matrix in $GL_2$ with $q^t=-q$. Then
   the Chern connection on $GL_2$ attached to $q$
    admits a correspondence structure $(\Gamma_p)$ with the following properties:  
   
   i) $\Gamma_p$ is irreducible, has  left degree $2$, and has right degree $2p^4$,
   
   ii) the  curvature of $(\Gamma_p)$ satisfies $\varphi^*_{pp'}=0$ for all $p,p'$, 
   
   iii)  the $(1,1)$-curvature of $(\Gamma_p)$ satisfies $\varphi^*_{p\overline{p}'}\neq 0$ for all $p,p'$.\end{theorem}
   
   \begin{theorem}
   \label{alta3}
    Let $q$ be the split matrix in $GL_2$ with $q^t=q$. Then
   the Chern connection on $GL_2$ attached to $q$
    admits a correspondence structure $(\Gamma_p)$ with the following properties:
   
   i)
   $\Gamma_p$ is irreducible and has left degree $4$, 
   
   ii)  the  $(1,1)$-curvature of $(\Gamma_p)$ satisfies $\varphi_{p\overline{p}'}^*\neq 0$ for all $p,p'$. 
 \end{theorem}

 Theorems \ref{alta2} and \ref{alta3} say that $Spec\ \bZ$ should be viewed as ``curved," in a natural way,  in the ``$(1,1)$-directions." 
 The vanishing of $\varphi_{pp'}^*$ Theorem \ref{alta2} should be viewed as a ``flatness" statement for the ``$(2,0)$ directions," in the antisymmetric case.
 Note that Theorem \ref{alta3}  says nothing about the curvature $\varphi_{pp'}^*$, i.e. about the ``$(2,0)$ directions,"  in the symmetric case; whether or not $\varphi_{pp'}^*$ vanishes in the case of Theorem \ref{alta3} is, at this point, an open problem and we have a similar open problem  in \cite{curvature1}.
 
The paper is organized as follows. 
Section 2 introduces some terminology related to correspondences on schemes. Section 3 specializes the discussion to correspondences on (spectra of) fields. Section 4 discusses compatibilities between correspondences on fields and lifts of Frobenius. Section 5 specializes the discussion to Chern connections and gives the proof  of our main results.

\bigskip

{\bf Acknowledgement}.
The author is indebted to J. Borger and Yu. I. Manin for inspiring conversations.
The  author would also like to acknowledge partial support from the Simons Foundation
(award 311773),  from the Institut des Hautes Etudes Scientifiques in Bures sur Yvette,   from the Romanian National Authority
for Scientific Research (CNCS - UEFISCDI, 
PN-II-ID-PCE-2012-4-0201), and from the Max-Planck-Institut f\"{u}r Mathematik in Bonn.

\section{Correspondences on a scheme}

In order to introduce our concepts  we need some terminology related to correspondences on schemes.  
The formalism below has a motivic flavor and can be viewed as a ``naive" variation on the Voevodsky formalism of finite correspondences \cite{voe}; this variation does not seem to naturally fit into the framework of \cite{voe} so our exposition will be independent of that in \cite{voe}.

A morphism of schemes will be called {\it totally dominant} if its image is dense and, moreover,  the image of any connected component of the source is dense in the corresponding connected component of the target; compositions of totally dominant morphisms are totally dominant.
Let $X$ be a  scheme.  By a {\it correspondence}
on $X$ we mean a triple $\Gamma=(Y,\pi,\varphi)$ where $\pi,\varphi:Y\ra X$ are morphisms of schemes. We sometimes write $Y=Y_{\Gamma}$; we also sometimes write
$\Gamma=\Gamma_{/X}=(Y_{/X},\pi_{/X},\varphi_{/X})$.
If $P$ is a property of schemes we say $\Gamma$ has $P$ if $Y_{\Gamma}$ has $P$. So we have a well defined notion, for instance, of {\it non-empty} correspondence, {\it connected} correspondence, {\it reduced} correspondence, etc. If $P$ is a property of morphisms of  schemes we say $\Gamma$ is {\it left} $P$ (respectively {\it right} $P$) if $\pi$ (respectively $\varphi$) has $P$.
So  we have a notion of {\it left \'{e}tale} correspondence, {\it right \'{e}tale} correspondence, {\it left totally dominant} correspondence, {\it right totally dominant} correspondence, {\it left finite} correspondence, {\it right finite} correspondence, etc. 
 There is a natural category ${\mathcal C}(X)$ of {\it correspondences} on $X$ whose objects are the correspondences on $X$; a morphism from $\Gamma'=(Y_{\Gamma'},\pi',\varphi')$
  to $\Gamma=(Y_{\Gamma},\pi,\varphi)$ is, by definition, a morphism $u:Y_{\Gamma'}\ra Y_{\Gamma}$ such that $\pi'=\pi\circ u$ and $\varphi'=\varphi\circ u$. In particular there is a notion of {\it isomorphism} of correspondences denoted by $\Gamma\simeq \Gamma'$. We say that a morphism of correspondences as above has property $P$  if $u$ has property $P$; in particular it is {\it totally dominant} if $u$ is totally dominant. 
   We say that a correspondence {\it totally dominates} another correspondence if there is a totally dominant morphism from the first to the second. 
  If $\Gamma=(Y,\pi,\varphi)$ is a correspondence we define its {\it transpose} as $\Gamma^t:=(Y,\varphi,\pi)$. If $\Gamma_1\ra \Gamma_2$ is a morphism of correspondences then we have an induced morphism $\Gamma_1^t\ra \Gamma_2^t$  and its formation is compatible with composition of morphisms.
  A correspondence $\Gamma=(Y,\pi,\varphi)$ will be called 
   {\it strictly symmetric} if $\Gamma^t=\Gamma$ i.e., if $\pi=\varphi$.
A correspondence    $\Gamma=(Y,\pi,\varphi)$ will be called 
  {\it symmetric} if $\Gamma^t\simeq \Gamma$, i.e., if there exists an isomorphism $\sigma:Y\ra Y$ such that $\varphi=\pi\circ \sigma$ and $\pi=\varphi\circ \sigma$.
   If $\Gamma\simeq\Gamma'$ and $\Gamma$ is symmetric  (respectively strictly symmetric) then $\Gamma'$ is also symmetric (respectively strictly symmetric). As an example, if  $Y \subset X \times X$ is a closed subscheme,  $\pi,\varphi:Y \ra X$ are given by the two projections from $X\times X$ to $X$, and $\tau$ is the automorphism of $X\times X$ that permutes the factors,  then the correspondence $(Y,\pi,\varphi)$ is symmetric if and only if $\tau(Y)=Y$; this correspondence is strictly symmetric if and only if $Y$ is contained in the diagonal.
  Going back to our general discussion, there are two remarkable correspondences on any $X$: the {\it empty} correspondence ${\mathbb O}:=(\emptyset,\emptyset,\emptyset)$ and the {\it identity} correspondence ${\mathbb I}=(X,id,id)$. A correspondence $(Y,\pi,\varphi)$  has a morphism to ${\mathbb I}$ if and only if it is strictly symmetric.  If $\Gamma=(Y_{\Gamma},\pi,\varphi)$ and $\Gamma'=(Y_{\Gamma'},\pi',\varphi')$
are two  correspondences on $X$ we define their composition as 
$$\Gamma'\circ \Gamma:=(Y_{\Gamma'\circ \Gamma},\pi\circ \pi'_1,\varphi'\circ \varphi_2)$$ where the above data are defined by the following diagram in which the square is cartesian:
$$\begin{array}{rclcc}
Y_{\Gamma'\circ \Gamma}
 & \stackrel{\varphi_2}{\ra} & Y_{\Gamma'} & \stackrel{\varphi'}{\ra} & X\\
\pi'_1 \downarrow & \  & \downarrow \pi' & \  & \  \\
Y_{\Gamma} & \stackrel{\varphi}{\ra} & X & \  & \  \\
\pi\downarrow & \  & \   & \ & \  \\
X & \  & \  & \  & \  \end{array}$$
Here, for any two morphisms of schemes with the same target we fix once and for all a fiber product object; the indices $1$ and $2$ stand for first and second projection.
In general $\Gamma'\circ \Gamma$ may be empty even if $\Gamma$ and $\Gamma'$ are non-empty; if, however, both correspondences are, say, left and right totally dominant and left and right of finite type  then their composition is non-empty.  Also, for two correspondences $\Gamma=(Y_{\Gamma},\pi,\varphi)$ and $\Gamma'=(Y_{\Gamma'},\pi',\varphi')$ we define their {\it direct sum} $\Gamma\oplus \Gamma'$ by taking
$Y_{\Gamma\oplus \Gamma'}:=Y_{\Gamma}\coprod Y_{\Gamma'}$
(disjoint union of schemes), with the obvious induced maps.
There are natural isomorphisms of correspondences expressing the ``associativity" of composition and direct sum, the ``neutral element" properties of ${\mathbb O}$ and ${\mathbb I}$, and the left and right ``distributivity" of composition with respect to direct sum. Cf. \cite{maclane}, pp. 158-159. All of these satisfy the standard (pentagonal and triangular) diagrams in loc. cit. plus diagrams for associativity, making  ${\mathcal C}(X)$  a monoidal  category  with respect to both composition and direct sum, having an additional  coherence structure for distributivity.
Also we have canonical isomorphisms 
$$(\Gamma\circ \Gamma')^t\simeq (\Gamma')^t\circ \Gamma^t,\ \ \ \Gamma\circ {\mathbb I}\simeq {\mathbb I}\circ \Gamma\simeq \Gamma,\ \ \ \Gamma\circ {\mathbb O}\simeq {\mathbb O}\circ \Gamma\simeq {\mathbb O},$$
satisfying appropriate coherence diagrams. 
If $\Gamma'\ra \Gamma''$ is a morphism of correspondences and $\Gamma$ is a correspondence we have natural morphisms of correspondences
\begin{equation}
\label{indu}
\Gamma'\oplus\Gamma\ra \Gamma''\oplus \Gamma,\ \ \ \Gamma\oplus \Gamma'\ra \Gamma\oplus \Gamma'',\ \ \ \Gamma'\circ \Gamma\ra \Gamma''\circ \Gamma,\ \ \ \Gamma\circ \Gamma'\ra \Gamma\circ \Gamma'',
\end{equation}
whose formation is compatible, in the obvious sense, with composition in ${\mathcal C}(X)$. The following definition will play a key role later:

\begin{definition}
\label{compatiblel}
A  correspondence $\Gamma$ is {\it  compatible} with a correspondence $\Gamma'$  if there is an isomorphism $\Gamma'\circ \Gamma^t\simeq \Gamma_1\oplus \Gamma_2$ where  $\Gamma_1$ is  strictly symmetric and left (equivalently, right) totally dominant. 
A correspondence $\Gamma_1$ is said to {\it partially commute}  with  $\Gamma_2$ if $\Gamma_1\circ \Gamma_2$ is compatible with   $\Gamma_2\circ \Gamma_1$; in other words, if the correspondence $\Gamma_2\circ \Gamma_1 \circ \Gamma_2^t \circ \Gamma_1^t$ has a direct summand which is strictly symmetric and left totally dominant. A correspondence $\Gamma_1$ is said to {\it  commute}  with  $\Gamma_2$ if $\Gamma_1\circ \Gamma_2$ is isomorphic to   $\Gamma_2\circ \Gamma_1$.
\end{definition}

If we consider the following diagram in which the square is cartesian,
\begin{equation}
\label{throat}\begin{array}{rclcc}
Y_{\Gamma'\circ \Gamma^t} & \stackrel{\pi_2}{\ra} & Y_{\Gamma'} & \stackrel{\varphi'}{\ra} & X\\
\pi'_1\downarrow & \  & \downarrow \pi' & \  & \  \\
Y_{\Gamma} & \stackrel{\pi}{\ra} & X & \  & \  \\
\varphi\downarrow & \  & \   & \ & \  \\
X & \  & \  & \  & \  \end{array}\end{equation}
the condition  that $\Gamma$ and $\Gamma'$ be   compatible is equivalent to  the condition that
 $\varphi'\circ \pi_2$ and $\varphi\circ \pi'_1$ coincide on {\it some} connected component  of $Y_{\Gamma'\circ \Gamma^t}$ and are totally dominant when restricted to that component.  Note that if $\Gamma$ is compatible with $\Gamma'$ then $\Gamma'$ is compatible  with $\Gamma$. Note that, in this generality, $\Gamma$ is not a priori compatible with itself; this is because $\Gamma\circ \Gamma^t$, although symmetric, does not a priori contain a strictly symmetric summand. However, if the image of the diagonal $Y_{\Gamma}\ra Y_{\Gamma\circ \Gamma^t}$ is open and closed then $\Gamma$ is compatible with itself; this happens, for instance, if $X$ is the spectrum of a field of characteristic zero and $\Gamma$ is reduced and left finite. 
 
 Similarly commutation does not imply partial commutation in general but it does imply it if $X$ is the spectrum of a field and we restrict ourselves to reduced,  left and right finite, correspondences.
 
  If $\sigma:X\ra X$ is an automorphism of $X$ we may consider the naturally associated correspondence $\Gamma_{\sigma}=(X,id,\sigma)$; we have $\Gamma_{\sigma}^t\simeq \Gamma_{\sigma^{-1}}$ and $\Gamma_{\sigma_1\circ \sigma_2}\simeq \Gamma_{\sigma_1}\circ \Gamma_{\sigma_2}$.

 If $X'\ra X$ is a morphism of schemes then we have a natural base change functor
 ${\mathcal C}(X)\ra {\mathcal C}(X')$, $\Gamma_{/X}\mapsto \Gamma_{/X'}$,
  compatible with all the operations introduced above. If $E$ is a field of characteristic zero we simply write ${\mathcal C}(E)$ in place of ${\mathcal C}(Spec\ E)$. 
 So, in particular, if $X$ is an  integral scheme of characteristic zero with field of  rational functions $E$ then $Spec\ E\ra X$ induces a functor 
 ${\mathcal C}(X)\ra {\mathcal C}(E)$, $\Gamma\mapsto \Gamma\otimes E$.

\section{Correspondences on a field}
Let $E$ be a field of characteristic zero.  
We denote by ${\mathcal C}_0(E)$ the subcategory of ${\mathcal C}(E)$ whose objects are the  non-empty, reduced,  left and right finite correspondences, and whose morphisms are the totally dominant (equivalently, surjective) morphisms. So the objects of ${\mathcal C}_0(E)$
 are of the form 
 $\Gamma=(Spec\ F,\pi,\varphi)$ where $F$ is a product of fields each of which is a finite extension of $E$ via  both $\pi$ and $\varphi$. We define the {\it bidegree} of $\Gamma$ as $\text{bideg}(\Gamma)=(\text{deg}(\pi),\text{deg}(\varphi))\in \bZ_{>0}\times \bZ_{>0}$; we then refer to $\text{deg}(\pi)$ and 
 $\text{deg}(\varphi)$ as the {\it left degree} and the {\it right degree} of $\Gamma$.
 One immediately checks that ${\mathcal C}_0(E)$ is closed under transposition $t$, composition $\circ$, taking direct sums $\oplus$, and taking direct summands.
 Moreover, if $\Gamma'\ra \Gamma''$  is a morphism in ${\mathcal C}_0(E)$ and $\Gamma$ is an object in ${\mathcal C}_0(E)$ then the morphisms in \ref{indu} are surjective, hence they are morphisms in ${\mathcal C}_0(E)$. Any morphism $\Gamma'\ra \Gamma$ in ${\mathcal C}_0(E)$ induces an injection $\cO(Y_{\Gamma})\ra \cO(Y_{\Gamma'})$; so  if there are morphisms $\Gamma'\ra \Gamma$ and $\Gamma\ra \Gamma'$ in 
 ${\mathcal C}_0(E)$ then $\Gamma\simeq \Gamma'$.  The relations  of compatibility  on the set of  objects in ${\mathcal C}_0(E)$ are reflexive and symmetric. As noted already, symmetry follows from the fact that, in this setting, $\Gamma\circ \Gamma^t$ contains a strictly symmetric summand namely, the image of the diagonal map $Y_{\Gamma}\ra Y_{\Gamma\circ \Gamma^t}$ which is both closed and open. Although compatibility is not transitive in general, there are cases when  transitivity holds. Indeed, assume that $\Gamma_i=(Y,\pi,\varphi_i)$, $i=1,2,3$, are correspondences in  ${\mathcal C}_0(E)$, with $Y$ irreducible, and $\pi$ Galois. It is easy to see that $\Gamma_i$ and $\Gamma_j$ are compatible if and only if there exists an automorphism $\sigma:Y\ra Y$ with $\pi\circ \sigma=\pi$ such that $\varphi_j=\varphi_i\circ \sigma$; hence if $\Gamma_1$ if compatible with $\Gamma_2$ and $\Gamma_2$ is compatible with $\Gamma_3$ then $\Gamma_1$ is compatible with $\Gamma_3$. 
 As already mentioned, if two correspondences in ${\mathcal C}_0(E)$ are commuting then they are partially commuting. Commutation of correspondences is, in principle a ``rare" occurrence; so it is important to put forward weaker concepts of commutation such as partial commutation (defined earlier) of the ``upper $*$" and ``lower $*$" commutation (to be defined presently).
 
 Let $\text{End}_{\text{gr}}(E)$ be the set of all group endomorphisms (equivalently of all ${\mathbb Q}$-linear endomorphisms) of   $E$; it is a ring with respect to addition and composition, and it is a ${\mathbb Q}$-Lie algebra in a natural way. (We used the subscript  $\text{gr}$ to avoid confusion with the monoid of field endomorphisms of $E$.)
 
 \begin{definition}
 Given a correspondence $\Gamma=(Y,\pi,\varphi)$ in ${\mathcal C}_0(E)$, $Y_{\Gamma}=Spec\ F$, and denoting by $\pi,\varphi:E\ra F$ the corresponding ring homomorphisms we define the map  $\Gamma^*\in \text{End}_{\text{gr}}(E)$, as the composition
 $\Gamma^*:=\text{tr}_{\pi}\circ \varphi: E\ra E \ra E$, 
 where $\text{tr}_{\pi}:F\ra E$ is the trace of $\pi$. \end{definition}
 
 It is clear that 
 $\Gamma^*=(\Gamma')^*$ for $\Gamma\simeq \Gamma'$ and 
 $(\Gamma_1\oplus \Gamma_2)^*=\Gamma_1^*+\Gamma^*_2$ in $\text{End}_{\text{gr}}(E)$; also  it is an easy exercise to check that $(\Gamma_1\circ \Gamma_2)^*=\Gamma_2^*\circ \Gamma^*_1$ in $\text{End}_{\text{gr}}(E)$. Note that if $\Gamma$ is strictly symmetric of bidegree $(m,m)$ then $\Gamma^*=m\cdot id$ with $id$ the identity on $E$.  
 More generally if $\Gamma'\ra \Gamma$ is a morphism  with $\Gamma$ irreducible and the morphism has degree $m$ then $(\Gamma')^*=m\cdot \Gamma^*$. Also for any $\Gamma$ of bidegree $(1,n)$, we have
 $(\Gamma^t)^*\Gamma^*=n\cdot id$ hence $\Gamma^*$   is injective. But if the left degree of $\Gamma$ is $\geq 2$ then $\Gamma^*$ is not injective, in general. Nevertheless, if $\Gamma$ is an arbitrary correspondence of left degree $m$ and $u\in E$ is such that $\Gamma^*(u^i)=0$ for all $1\leq i\leq m$ then $u=0$.
 
Assume (as will be the case in applications) that $E$ comes equipped with an involution $\iota:E\ra E$, i.e. a field automorphism whose square is the identity; we assume $\iota$ is not the identity. 
Consider the eigenspaces of $\iota$:
\begin{equation}
\label{E+}
E^+=\{z\in E;\iota(z)=z\},\ \ E^-=\{z\in E;\iota(z)=- z\};\end{equation}
then $E$ has degree $2$ over the field $E^+$ and we have an $E^+$-linear space decomposition $E=E^+\oplus E^-$. Say that a correspondence $\Gamma$ on $E$ is {\it partially induced} from $E^+$ if $\Gamma^*(E^+)\subset E^+$ and $\Gamma^*(E^-)=0$.
If, in addition, $\varphi(E^+)\subset \pi(E^+)$ we say that $\Gamma$ is {\it induced} from $E^+$; if this is the case then, of course, the restriction of $\Gamma^*$ to $E^+$ is an integer multiple of a field endomorphism of $E^+$ and, in particular, $Ker(\Gamma^*)=E^-$.

\begin{remark}
The action of correspondences on fields extends to an action on forms as follows.
Assume, for simplicity, that $E$ is a finitely generated field (as usual, of characteristic zero). For any finite reduced $E$-algebra $\pi:E\ra F$ let $\Omega_F$ be the $F$-module of K\"{a}hler differentials of $F$ over ${\mathbb Q}$ and $\Omega_F^i$ its $i$-th exterior power. We have natural isomorphisms $\Omega^i_E\otimes F\simeq \Omega^i_F$
and hence naturally induced maps 
$$\text{tr}^i_{\pi}:\Omega^i_F\simeq \Omega^i_E\otimes_E F\stackrel{1\otimes tr_{\pi}}{\longrightarrow} \Omega^i_E\otimes_E E\simeq \Omega^i_E.$$
If $x_1,...,x_d$ is a transcendence basis of $E$ over ${\mathbb Q}$ and $\frac{\partial}{\partial x_i}$ are the corresponding derivations on $E$ and $F$ then, 
$\Omega_F$ is a free $F$-module with basis $dx_1,...,dx_d$ and,
for any $f\in F$, $df=\sum \frac{\partial f}{\partial x_i}dx_i$, hence
$$\text{tr}^1_{\pi}(df)=\sum \text{tr}_{\pi}\left(\frac{\partial f}{\partial x_i}\right)dx_i=\sum \frac{\partial}{\partial x_i}(\text{tr}_{\pi} f)dx_i=d(\text{tr}_{\pi} f).$$
(The commutation of $\text{tr}_{\pi}$ and $\frac{\partial}{\partial x_i}$ is checked by reducing to the case when $F$ is a Galois field extension of $E$ in which case it follows from the commutation of $\frac{\partial}{\partial x_i}$ with the $E$-automorphisms of $F$.) Now if 
$\text{End}(\Omega_E)$ denotes the ring of group endomorphisms (equivalently ${\mathbb Q}$-linear endomorphisms) of $\Omega_E$ and 
$\Gamma=(Spec\ F,\pi,\varphi)$ is a correspondence in ${\mathcal C}_0(E)$ we can define a ${\mathbb Q}$-linear map $\Gamma^{*i}\in \text{End}(\Omega_E)$, $\Gamma^{*i}:\Omega^i_E\ra \Omega^i_E$ by the formula $\Gamma^{*i}=\text{tr}^i_{\pi}\circ \varphi^{*i}:\Omega^i_E\ra \Omega^i_F\ra \Omega^i_E$, where $\varphi^{*i}$ is naturally induced from $\varphi$. One has $\Gamma^{*1}(df)=d(\Gamma^*(f))$ for all $f\in F$. Also for any two correspondences $\Gamma_1,\Gamma_2$ we have 
$(\Gamma_1\oplus \Gamma_2)^{*i}=\Gamma^{*i}_1\oplus \Gamma^{*i}_2$ and 
$(\Gamma_1\circ \Gamma_2)^{*i}=\Gamma^{*i}_2\circ \Gamma^{*i}_1$. So if $\Gamma_1$ and $\Gamma_2$ commute then $\Gamma_1^{*i}$ and $\Gamma^{*i}_2$ commute.

If we assume in addition that an involution $\iota:E\ra E$ is given such that $\iota(x_i)=-x_i$ (which is the case in our applications) and if $\Gamma$ is partially induced from $E^+$ (in particular $\Gamma^*(E^-)=0$) then we have $\Gamma^{*1}(dx_i)=d(\Gamma^*(x_i))=0$. 
\end{remark}

  There is an alternative action of correspondences on ``cycles" of $E$ as follows.
   Fix an algebraically closed field  ${\mathbb U}$, of cardinality at least that of $E$, let $\text{Hom}_{\text{fields}}(E,{\mathbb U})$ denote
   the set all field homomorphisms $\sigma:E\ra {\mathbb U}$, and denote by $\text{Cyc}(E)$ the free abelian group with basis $\text{Hom}_{\text{fields}}(E,{\mathbb U})$. 
   For $\sigma\in \text{Hom}_{\text{fields}}(E,{\mathbb U})$ denote by $[\sigma]$ the image of $\sigma$ in $\text{Cyc}(E)$. The inclusion $\text{Hom}_{\text{fields}}(E,{\mathbb U})\ra \text{Hom}_{\text{gr}}(E,{\mathbb U})$ induces a group homomorphism $\int:\text{Cyc}(E)\ra \text{Hom}_{\text{gr}}(E,{\mathbb U})$,
   $$\int\left(\sum n_i[\sigma_i]\right):=\sum n_i \sigma_i.$$
    Let $\text{End}(\text{Cyc}(E))$ be the ring of group endomorphisms of $\text{Cyc}(E)$. 
   For any correspondence $\Gamma$ in ${\mathcal C}_0(E)$ we can define an endomorphism
   $\Gamma_*\in \text{End}(\text{Cyc}(E))$ as follows.
   
   \begin{definition}
    Given a correspondence $\Gamma$ in ${\mathcal C}_0(E)$ we define the map  $\Gamma^*\in \text{End}_{\text{gr}}(\text{Cyc}(E))$ as follows. Assume first
   $\Gamma=(Y,\pi,\varphi)$ is irreducible, so $Y_{\Gamma}=Spec\ F$ where $F$ is a field and $\pi,\varphi$ are induced by finite field homomorphisms (still denoted by) $\pi,\varphi:E\ra F$. Let $\sigma:E\ra {\mathbb U}$ be an embedding. Then consider all the embeddings $\sigma_1,...,\sigma_d:F\ra {\mathbb U}$ such that
 $\sigma_i\circ \pi=\sigma$, where $d$ is the degree of $\pi$, and define
 $$\Gamma_*[\sigma]=\sum_{i=1}^d [\sigma_i\circ \varphi]\in \text{Cyc}(E).$$
 Extend $\Gamma_*$ by linearity to an endomorphism $\Gamma_*$ of $\text{Cyc}(E)$. Finally extend the definition of $\Gamma_*$ by linearity to the case when $\Gamma$ is not necessarily irreducible. \end{definition}
 
 Clearly, if $\Gamma\simeq \Gamma'$ then $\Gamma_*=\Gamma'_*$.  By definition $(\Gamma\oplus \Gamma')_*=\Gamma_*+\Gamma'_*$ and  it is a trivial exercise to check that
 $(\Gamma\circ \Gamma')_*=\Gamma_*\circ \Gamma'_*$. 
 So if $\Gamma_1$ and $\Gamma_2$ commute then $\Gamma_{1*}$ and $\Gamma_{2*}$ commute.
 Also it is trivial to check that if $\Gamma$ is strictly symmetric of bidegree $(d,d)$ then
  $\Gamma_*=d$, where $d$ is the multiplication by $d$ endomorphism of $\text{Cyc}(E)$. 
  
  \medskip

The above construction $\Gamma\mapsto \Gamma_*\in \text{End}(\text{Cyc}(E))$ is related to the previous construction $\Gamma\mapsto \Gamma^*\in \text{End}_{\text{gr}}(E)\subset \text{Hom}_{\text{gr}}(E,{\mathbb U})$ as follows: for any correspondence $\Gamma$ on $E$ and any $\sigma\in \text{Hom}_{\text{fields}}(E,{\mathbb U})$ we have an equality
\begin{equation}
\label{compatty}
\int\left(\Gamma_* [\sigma]\right)=\sigma\circ \Gamma^*
\end{equation}
in $\text{Hom}_{\text{gr}}(E,{\mathbb U})$. In particular, if $\Gamma_1,\Gamma_2$ are two correspondences and 
$\Gamma_{1*}= \Gamma_{2*}$ then 
$\Gamma_1^*=\Gamma_2^*$.  So if $\Gamma,\Gamma'$ are two correspondences 
such that $\Gamma_*$ and $\Gamma'_*$ commute, i.e., 
$[\Gamma_*,\Gamma'_*]= 0$ in $\text{End}(\text{Cyc}(E))$ 
then $\Gamma^*$ and $(\Gamma')^*$ commute, i.e., 
$[\Gamma^*,(\Gamma')^*]= 0$ in $\text{End}_{\text{gr}}(E)$.  
 We summarize the various implications between various commutation relations:
 $$\begin{array}{ccc}
 \text{$\Gamma$ and $\Gamma'$ commute} &  \Longrightarrow & \text{$\Gamma_*$ and $\Gamma'_*$ commute}\\
 \  & \  & \  \\
 \Downarrow & \  & \Downarrow\\
 \  & \  & \  \\
  \text{$\Gamma$ and $\Gamma'$ partially commute} &  \  & \text{$\Gamma^*$ and $(\Gamma')^*$ commute}\end{array}
 $$

\begin{remark}
Here is a ring theoretic formalism into which the above definitions fit; this will not be needed in the sequel but seems to provide the right context for our concepts.
Assume, as before, that $E$ is a field of characteristic zero. We can associate to $E$ a ring $C=C(E)$ equipped with  an anti-involution as follows. Let us denote by $C_+=C_+(E)$ the set of isomorphism classes $[\Gamma]$ of objects $\Gamma$  of ${\mathcal C}_0(E)$
  to which we add one element $[{\mathbb O}]$ which we call {\it zero} and we denote by $0$. We denote the element $[{\mathbb I}]$ by $1$.
 The set $C_+$ comes equipped with the following operations:
 1) transposition: $[\Gamma]^t:=[\Gamma^t]$;
 2)  addition: $[\Gamma]+ [\Gamma']:=[\Gamma\oplus\Gamma']$;
 3) multiplication: $[\Gamma]\cdot[\Gamma']:=[\Gamma\circ\Gamma']$.
 These operations make $C_+$ into a {\it semiring with anti-involution} by which we mean that the following hold:
 i) addition and multiplication are associative and  addition is commutative;
 ii) multiplication is left and right distributive;
 iii) transposition is an anti-involution (i.e. it is a homomorphism with respect to direct sum, an anti-homomorphism with respect to composition, and it is its own inverse);  
 iv) $0$ is a neutral element for addition, $0^t=0$,  and  the product of $0$ with any element is $0$; 
  v)  
  $1$ is a neutral element for composition and $1^t=1$. In addition the sum and product of two non-zero elements in $C_+$
 is non-zero.
  Let us say that an element of $C_+$ is {\it irreducible} if it is non-zero and cannot be written as a sum of two non-zero elements.
 Clearly $[\Gamma]$ is irreducible if and only if $Y_{\Gamma}$ is non-empty and irreducible. Also any non-zero element of $C_+$ can be written as a finite sum of irreducible elements and this decomposition is unique up to the permutation of the terms. In other words $C_+$ with its addition is a free commutative monoid with basis the irreducible elements. In particular $C_+$ has the  cancellation property for addition. We can now embed our semiring above into a ring as follows. Define
 $C(E):=C:=(C_+\times C_+)/\sim$
 where $(\gamma_1,\gamma_2)\sim (\gamma_3,\gamma_4)$ if and only if
 $\gamma_1+\gamma_4=\gamma_2+\gamma_3$. Due to the cancellation property $\sim$ is an equivalence relation and the map 
 $C_+\ra C$, $\gamma\mapsto (\gamma,0)$,
 is injective. Then $C$ has a naturally induced structure of ordered ring with identity and is equipped with an induced  anti-involution.  Also $C$ is a free $\bZ$-module with respect to addition, with basis the irreducible elements of $C_+$. 
The ring $C$ can be referred to as the {\it ring of correspondences} on $E$.
Let $\text{End}_{\text{field}}(E)$ be the monoid of field endomorphisms of $E$ and $\text{Aut}_{\text{field}}(E)$ be the group of field automorphisms of $E$.
 Note  that one has a natural injective homomorphisms of multiplicative monoids 
 $\text{End}_{\text{field}}(E)\ra C(E)$ given by $\sigma\mapsto [\Gamma_{\sigma}]$, $\Gamma_{\sigma}=(Spec\ E,id,\sigma)$; if one views this embedding as an inclusion then, for any $$\sigma\in \text{Aut}_{\text{field}}(E)\subset \text{End}_{\text{field}}(E)\subset C(E)$$
  we have $\sigma^t=\sigma^{-1}$. 
  One has a unique ring (anti)homomorphism 
  $$(\ )^*:C(E)\ra \text{End}_{\text{gr}}(E)$$ such that $[\Gamma]\mapsto \Gamma^*$ for any correspondence $\Gamma$. More generally one has induced ring homomorphisms $C(E)\ra \text{End}(\Omega^i_E)$ sending each $[\Gamma]$ into $\Gamma^{*i}$.
  The commutator $[\Gamma^*_1,\Gamma^*_2]$ in $\text{End}_{\text{gr}}(E)$, that will soon play a role in our theory of curvature, appears then as the image in $\text{End}_{\text{gr}}(E)$ of negative of the commutator
  $[\Gamma_1][\Gamma_2]-[\Gamma_2][\Gamma_1]$ in $C(E)$; and the same is the case for the commutators 
  $[\Gamma^{*i}_1,\Gamma^{*i}_2]$ in $\text{End}(\Omega_E^i)$.
  On the other hand there is a natural ring homomorphism 
  $$(\ )_*:C(E)\ra \text{End}(\text{Cyc}(E))$$ sending any $[\Gamma]$ into $\Gamma_*$ hence the commutator 
  $[\Gamma_{1*},\Gamma_{2*}]$ in $\text{End}_{\text{gr}}(\text{Cyc}(E))$ is the image  of the commutator
  $[\Gamma_1][\Gamma_2]-[\Gamma_2][\Gamma_1]$ in $C(E)$.
     Note that by \ref{compatty} we have, for any embedding $\sigma:E\ra {\mathbb U}$, a commutative diagram of groups
   $$\begin{array}{rcl}
   C(E) & \stackrel{(\ )^*}{\longrightarrow} & \text{End}_{\text{gr}}(E)\\
   \  & \  & \  \\
   (\ )_* \downarrow & \  & \downarrow \sigma\circ  \\
   \  & \  & \  \\
   \text{End}(\text{Cyc}(E)) & \stackrel{\int((\ )[\sigma])}{\longrightarrow} & \text{Hom}_{\text{gr}}(E,{\mathbb U}).\end{array}
   $$
   We end by noting that it would be interesting to answer the following questions. Is the subset of $C(E)$ consisting of all elements of the form $[\Gamma]$ with $\Gamma$ strictly symmetric  definable in $C(E)$ in terms of the ring operations and the transposition operation only?  Is partial commutation in $C(E)$  definable in terms of these operations only? What are the kernels of the ring homomorphisms $(\ )^*$ and  $(\ )_*$? Since the right vertical homomorphism in the above diagram is injective we get, of course, that $Ker((\ )_*)\subset Ker((\ )^*)$.
   \end{remark}

\section{Correspondences compatible with lifts of Frobenius}
Let $B$ be a ring and $p$ a prime integer. By a {\it lift of Frobenius} on $B$ we mean a ring endomorphism 
of $B$ whose reduction mod $p$ is the $p$-power Frobenius endomorphism of $B/pB$.
We denote by $B^{\widehat{p}}$ the $p$-adic completion of $B$. For $X$ a noetherian scheme we denote by $X^{\widehat{p}}$ the $p$-adic completion of $X$, which is a $p$-adic formal scheme. By a {\it lift of Frobenius} on $X$, respectively $X^{\widehat{p}}$, we mean an endomorphism of $X$ or $X^{\widehat{p}}$ whose reduction mod $p$ is the $p$-power Frobenius.

Fix, in what follows, a collection ${\mathcal V}$ of prime integers.
 
Now let us assume $X$ is an affine smooth connected scheme over a regular integral domain $A$. Note that $X$ is then integral; we denote by $K$ the fraction field of $A$ and we denote by $E$ the fraction field of $X$. (These hypotheses  will always be satisfied in our applications which deal with the case $X=GL_n$ over $A$.) 

\begin{definition}
\label{soso}
We say that a lift of Frobenius $\phi_p$ on $X^{\widehat{p}}$ and 
 a  correspondence $\Gamma_{p}=(Y_{p},\pi_{p},\varphi_{p})$ in ${\mathcal C}(E)$ are {\it  compatible} 
if   there exists a correspondence 
$\Gamma_{p/X}=(Y_{p/X},\pi_{p/X},\varphi_{p/X})$ on $X$ such that the following properties are satisfied.

1) $\Gamma_{p}\simeq \Gamma_{p/X}\otimes E$;

2) $\Gamma_{p/X}$ is affine and left \'{e}tale;

3) $\pi_{p/X}^{\widehat{p}}:Y_{p/X}^{\widehat{p}}  \ra X^{\widehat{p}}$ is an isomorphism;

4) $\varphi_{p/X}^{\widehat{p}}=
\phi_p\circ \pi_{p/X}^{\widehat{p}}:Y_{p/X}^{\widehat{p}}\ra X^{\widehat{p}}$.
\end{definition}

If $\phi_p$ is compatible with a correspondence $\Gamma_{p}$ then the correspondence is, of course, not unique; but it is ``essentially unique" in the following sense. 

\begin{lemma}
\label{abovel}
Assume  that $\phi_p$ is  compatible with a  correspondence $\Gamma_{p}$ in ${\mathcal C}(E)$. Then $\Gamma_{p}$ is in ${\mathcal C}_0(E)$.  
Assume   in addition that
$X\otimes {\mathbb F}_p$ is connected for all $p\in \cP$ and  $\phi_p$ is also  compatible with a correspondence $\Gamma'_{p}$ in ${\mathcal C}(E)$; then
$\Gamma_{p}$ and $\Gamma'_{p}$ are  compatible. 
\end{lemma}

\begin{remark}
\label{abovell}
In our applications 
$$A=\bZ[1/M,\zeta_N],\ \ \ X=GL_{n/A}:=Spec\ A[x,\det(x)^{-1}].$$
 So, for $p\not| MN$,  the condition that $X\otimes {\mathbb F}_p$ be connected 
is equivalent to $p$ being inert in ${\mathbb Q}(\zeta_N)$. \end{remark}

{\it Proof of Lemma \ref{abovel}}.
Let  $\Gamma_{p/X}=(Y_{p/X},\pi_{p/X},\varphi_{p/X})$ be as in Definition \ref{soso} and let
$Y$ be a connected component of $Y_{p/X}$ that meets $Y_{p/X}\otimes {\mathbb F}_p$. Then
$Y$ is smooth over $A$, hence $Y$ is a regular scheme, hence an integral scheme (because it is connected), hence    $\cO(Y)$ is $p$-adically separated (by Krull's intersection theorem, because $p$ is not invertible on $Y$) and therefore $\cO(Y)$ embeds into $\cO(Y)^{\widehat{p}}$. But  the maps $\pi_{p/X}^{\widehat{p}},\varphi_{p/X}^{\widehat{p}}:\cO(X)^{\widehat{p}}\ra \cO(Y)^{\widehat{p}}$ are injective (because their reductions mod $p$ are injective and the rings in question are $p$-adically complete with $p$ a non-zero divisor in them). So the maps
$\pi_{p/X},\varphi_{p/X}:\cO(X)\ra \cO(Y)$ are injective. This immediately implies that $\Gamma_{p}$ is non-empty. Since $\pi_{p/X}$ is \'{e}tale it is generically finite so
$\cO(Y_{p})$ is a product of fields which are finite extensions of $E$. Due to considerations of transcendence degrees we get that $\varphi_{p}$ is finite. Hence  
$\Gamma_{p}$ is non-empty, reduced, and  left and right finite.

Assume now $\phi_p$ is also  compatible with a correspondence 
$\Gamma'_{p}$ and let 
  $\Gamma'_{p/X}=(Y'_{p/X},\pi'_{p/X},\varphi'_{p/X})$ be as in Definition \ref{soso}. 
Let  $Y$  be a connected component of $Y_{\Gamma'_{p/X}\circ \Gamma_{p/X}^t}$ that meets $Y_{\Gamma_{p/X}'\circ \Gamma_{p/X}^t}\otimes {\mathbb F}_p$ (the latter being non-empty because $X\otimes {\mathbb F}_p$ is irreducible). Then, as before,
$\cO(X)$ and $\cO(Y)$ embed  into their $p$-adic completions; on the other hand, with notation as in diagram \ref{throat} (with $\pi=\pi_{p/X}$,
$\pi'=\pi'_{p/X}$, $\varphi=\varphi_{p/X}$, $\varphi'=\varphi'_{p/X}$)  the $p$-adic completions of $\varphi'\circ \pi_2$ and $\varphi\circ \pi'_1$  are equal. So $\varphi'\circ \pi_2$ and $\varphi\circ \pi'_1$  are equal on $Y$ hence $\Gamma_{p/X}$ and $\Gamma'_{p/X}$ are  compatible. Since, as before, $\pi'\circ \pi_2:\cO(X)\ra \cO(Y)$ is injective it follows that $Y\otimes E$ is non-empty, hence
$\Gamma_{p}$ and $\Gamma'_{p}$ are  compatible.
\qed

\medskip

In what follows we specialize our discussion to $X=GL_n$. Recall first that a $p$-derivation on a ring $B$ is a map $\d_p:B\ra B$ such that the map $\phi_p:B\ra B$, $\phi_p(b)=b^p+p\d_p(b)$ is a ring homomorphism (referred to as the lift of Frobenius attached to $\d_p$); cf. the Introduction. Then we fix, for the rest of the paper,  the following terminology and notation:

\begin{definition}
\label{basique}
Let $A=\bZ[1/M,\zeta_N]$, with $M$ even, set  
$$G=GL_n=Spec\ A[x,\det(x)^{-1}],\ \ x=(x_{kl})_{1\leq k,l\leq n},$$  let $K={\mathbb Q}(\zeta_N)$ denote the fraction field of $A$ and we let $E=K(x)$ denote the fraction field of $GL_n$.
Let ${\mathcal V}$ be a set of primes not dividing $MN$.
An adelic connection on $GL_n$ is a family $\d=(\d_p)$ indexed by ${\mathcal V}$ where, for each $p$, $\d_p$ is a $p$-derivation on $A[x\det(x)^{-1}]^{\widehat{p}}$. 
 Let  $\d=(\d_p)$ be an adelic connection on $GL_n$ 
 with associated family $(\phi_p)$ of lifts of Frobenius.
 By a {\it correspondence structure} for $\d=(\d_p)$ 
we understand a family $(\Gamma_{p})$ of correspondences in ${\mathcal C}_0(E)$ such that, for each $p$, $\phi_p$ is  compatible with $\Gamma_p$.  We define the 
(``upper $*$") 
{\it curvature} of $(\Gamma_p)$ as 
the matrix  $(\varphi^*_{pp'})$, $\varphi^*_{pp'}=\frac{1}{pp'}[\Gamma_{p'}^*,\Gamma^*_p]\in \text{End}_{\text{gr}}(E)$; cf. \ref{varfistar1}. Given one more adelic connection 
$(\d_{\overline{p}})$ with correspondence structure $(\Gamma_{\overline{p}})$ we define the 
(``upper $*$") 
$(1,1)$-{\it curvature} of $(\Gamma_p)$  with respect to $(\Gamma_{\overline{p}})$ as the matrix $(\varphi^*_{p\overline{p}'})$, $\varphi^*_{p\overline{p}'}=\frac{1}{pp'}[\Gamma^*_{\overline{p}'},\Gamma^*_p]\in \text{End}_{\text{gr}}(E)$ for $p\neq p'$ and 
$\varphi^*_{p\overline{p}}=\frac{1}{p}[\Gamma^*_{\overline{p}},\Gamma^*_p]\in \text{End}_{\text{gr}}(E)$; cf. \ref{varfistar2}, \ref{varfistar22}.
Similarly we define the (``lower $*$") curvature and $(1,1)$-curvature $(\varphi_{pp'*})$, $\varphi_{pp'*}=\frac{1}{pp'}[\Gamma_{p'*},\Gamma_{p*}]\in \text{End}_{\text{gr}}(\text{Cyc}(E))$ and $(\varphi_{p\overline{p}'*})$, $\varphi_{p\overline{p}'*}=\frac{1}{pp'}[\Gamma_{\overline{p}'*},\Gamma_{p*}]\in \text{End}_{\text{gr}}(\text{Cyc}(E))$ for $p\neq p'$ and 
$\varphi_{p\overline{p}*}=\frac{1}{p}[\Gamma_{\overline{p}*},\Gamma_{p*}]\in \text{End}_{\text{gr}}(\text{Cyc}(E))$. (Note that, by \ref{compatty}, the non-vanishing of the upper $*$ curvature, respectively $(1,1)$-curvature, implies the non-vanishing of the corresponding lower $*$  curvature.)
\end{definition}

\begin{remark}
Correspondence structures  are not unique; but, in our applications  they are ``essentially unique" in the sense of Lemma \ref{abovel}.\end{remark}

\begin{example}
\label{caldura}
Consider the adelic connection $\d_0=(\d_{0p})=\overline{\d}=(\d_{\overline{p}})$ on $G=GL_n$ with $\d_{0p}x=\d_{\overline{p}}x=0$, indexed by a set of odd primes. The attached family of lifts of Frobenius $(\phi_{0p})$
is given by $\phi_{0p}(x)=x^{(p)}:=(x_{ij}^p)$. Then $\d_0=\overline{\d}$
 admits a   correspondence structure which we denote by $(\Gamma_{\overline{p}})$, referred to as the {\it canonical} correspondence structure, as follows: take  $\Gamma_{\overline{p}}=\Gamma_{\overline{p}/G}\otimes E$, where $\Gamma_{\overline{p}/G}=(Y_{\overline{p}/G},\pi_{\overline{p}/G},\varphi_{\overline{p}/G})$, 
$$Y_{\overline{p}/G}=Spec\ A[x,\det(x)^{-1},\det(x^{(p)})^{-1}],$$
 $\pi_{\overline{p}/G}$ is induced by the natural inclusion $\cO(G)\subset \cO(Y_p)$, and $\varphi_{\overline{p}/G}$
is induced by the map $\cO(G)\ra\cO(Y_p)$, $x\mapsto x^{(p)}$. Clearly $\Gamma_{\overline{p}}$ and $\Gamma_{\overline{p}'}$ commute for all $p,p'$. In particular the  ``lower $*$" and ``upper $*$" curvatures of $\overline{\d}$ vanish.
\end{example}

\begin{definition}
 Let $\iota: E\ra E=K(x)$ be the $K$-automorphism with $\iota(x)=- x$. 
and recall that we denote by $E^+$ and $E^-$ the eigenspaces of $\iota$ in $E$.
A correspondence $\Gamma$ is said to be {\it partially induced} from $E^+$  if $\Gamma^*(E^+)\subset E^+$ and $\Gamma^*(E^-)=0$; it is called
 {\it induced} from $E^+$
if it is partially induced from $E^+$ and $\varphi(E^+)\subset \pi(E^+)$.\end{definition}

\begin{remark}
In the definition above $E$ has degree $2$ over its subfield $E^+$,  $E^+$ is generated as a field over $K$ by the monomials of degree $2$ in the entries of $x$, and  $E=E^+\oplus E^-$. 
\end{remark}

\begin{remark}
For the canonical correspondence structure attached to $(\d_{\overline{p}})$, we have $\Gamma_{\overline{p}}^*(E^+)\subset E^+$ and $\Gamma^*_{\overline{p}}(E^-)\subset E^-$. So if a correspondence structure $(\Gamma_p)$ on an adelic connection $(\d_p)$ has the property that all $\Gamma_p$'s are partially induced from $E^+$ then the curvature and $(1,1)$-curvature have the property that 
$$\varphi^*_{pp'}(E^+)\subset E^+,\ \ \varphi^*_{p\overline{p}'}(E^+)\subset E^+,\ \ \varphi^*_{pp'}(E^-)=0,\ \ \varphi^*_{p\overline{p}'}(E^-)=0.$$
In particular the curvature and $(1,1)$-curvature are completely determined by their restriction to $E^+$. \end{remark}

\begin{remark} In the notation of Definition \ref{basique} one can define, more generally, 
matrices  $(\varphi^{*i}_{pp'})$, $\varphi^{*i}_{pp'}=\frac{1}{pp'}[\Gamma_{p'}^{*i},\Gamma^{*i}_p]\in \text{End}(\Omega^i_E)$ and $(\varphi^{*i}_{p\overline{p}'})$, $\varphi^{*i}_{p\overline{p}'}=\frac{1}{pp'}[\Gamma^{*i}_{\overline{p}'},\Gamma^{*i}_p]\in \text{End}(\Omega^i_E)$ for $p\neq p'$ and 
$\varphi^{*i}_{p\overline{p}}=\frac{1}{p}[\Gamma^{*i}_{\overline{p}},\Gamma^{*i}_p]\in \text{End}(\Omega^i_E)$. Then, by analogy with \cite{curvature1}, one can try to use these matrices (for $i=1$) to define {\it first Chern forms} and {\it first Chern $(1,1)$-forms} respectively. One can consider, for instance,  the elements
$\varphi^{*1}_{pp'}(dx_{ij})=\sum_{kl} a_{ijkl,pp'}dx_{kl}$ and define, say, the first Chern form as the matrix $(\rho_{pp'})$ with entries $\rho_{pp'}=\sum_{ij}a_{ijij,pp'}$. 
With this definition, however, if $\Gamma_p$ are partially induced from $E^+$ (which will be the case in applications, cf. Theorem \ref{rational} below), we have $\Gamma^{*1}_p(dx_{ij})=0$ for all $p,i,j$ so we get $\rho_{pp'}=0$ for all $p,p'$.
One can also consider the elements
$\varphi^{*1}_{p\overline{p}}(dx_{ij})=\sum_{kl} a_{ijkl,p\overline{p}}dx_{kl}$ and define, say, the $(p,\overline{p})$- components of a ``first Chern $(1,1)$-form" by $\rho_{p\overline{p}}=\frac{1}{p}\sum_{ij}a_{ijij,p\overline{p}}$. Since $\Gamma^{*1}_{\overline{p}}(dx_{ij})=d(x_{ij}^p)$ and $x_{ij}^p\in E^-$ it follows again that
$\rho_{p\overline{p}}=0$. The above comments show that a natural definition of first Chern forms should, maybe, take into account the action of our correspondences on forms in $dE^+$ rather than on the forms $dx_{ij}\in dE^-$. We will not pursue the discussion on first Chern forms further in this paper.
\end{remark}

\section{Correspondence structure for Chern connections}
We assume the notation of Definition \ref{basique} and Example \ref{caldura}. 
Recall from the Introduction the following: 

\begin{definition}
 Let $q\in GL_n(A)$ with $q^t=\pm q$.  The {\it Chern connection} attached to $q$ is the unique adelic connection $\d=(\d_p)$ on $GL_n$ such that if $(\phi_p)$ is the attached family of lifts of Frobenius then $\phi_p$ and $\phi_{0p}$  make the diagrams \ref{coconut} commutative; $\d$ exists and is unique by \cite{adel2}.\end{definition}

We recall from \cite{adel2,curvature1} that 
if $\d=(\d_p)$ is the Chern connection attached to $q$ and $(\phi_p)$ are 
the lifts of Frobenius  on $\cO(G)^{\widehat{p}}=A[x\det(x)^{-1}]^{\widehat{p}}$ attached to  $\d$ then for each $p$, $\phi_p$ is  the unique    ring homomorphism
$$\phi_p:\cO(G)^{\widehat{p}}\ra \cO(G)^{\widehat{p}},\ \ \ x\mapsto \phi_p(x)=\Phi_p,$$
such that $\phi_p(\zeta_N)=\zeta_N^p$ and
 \begin{equation}
 \label{turandot}
\Phi_p(x)=x^{(p)}\{(x^{(p)t}\phi_p(q)x^{(p)})^{-1} (x^tqx)^{(p)}\}^{1/2}.\end{equation}
In this expression, the  exponent $t$ for a matrix means {\it transpose},  the upper index $(p)$ for a matrix  denotes the operation of  raising all the entries of the matrix to the $p$-th power,  and the $1/2$ power is computed using the usual matrix series  
\begin{equation}
\label{radical}
(1+u)^{1/2}=1+\sum_{i=1}^{\infty} \left(\begin{array}{c} 1/2\\i\end{array}\right) u^i\end{equation}
 for  $u$ an $n\times n$ matrix of indeterminates. The formula is applied to
 $$u=(x^{(p)t}\phi_p(q)x^{(p)})^{-1} (x^tqx)^{(p)}-1,$$
 which has coefficients in $p\cO(G)$; this makes the series \ref{radical} convergent because
   $\left(\begin{array}{c} 1/2\\i\end{array}\right)\in \bZ[1/2]$.
  
 Let us say that a matrix $q\in GL_{2r}(\bZ)$ is  {\it split} if it is one of the following:
 \begin{equation}
 \label{scorpion3}
  \left(\begin{array}{cl} 0 & 1_r \\-1_r & 0\end{array}\right),\ \ 
\left( 
\begin{array}{ll} 0 & 1_r\\1_r & 0\end{array}\right),
\end{equation}
where $1_r$ is the $r \times r$ identity matrix.

\begin{theorem}
\label{rational}\ 

1)
For any  $q\in GL_n(A)$ with $q^t=\pm q$ the Chern connection on $GL_n$ attached to $q$ admits a   correspondence structure.

2) If  $q\in GL_2(A)$ is split, with $q^t=-q$, then the Chern connection on $GL_2$ attached to $q$ admits a correspondence structure  $(\Gamma_p)$ with the following properties:
  each $\Gamma_p$ is irreducible, is induced from $E^+$, and  has left degree $2$ and right degree $2p^4$;  $(\Gamma_p)$ has curvature satisfying $\varphi_{pp'}^*=0$ and $(1,1)$-curvature satisfying $\varphi^*_{p\overline{p}'}\neq0$,  hence  $\varphi_{p\overline{p}'*}\neq0$,   for all $p,p'$; and for all $p,p'$,  $\Gamma_p$ and $\Gamma_{p'}$ partially commute, but do not commute.

3) If $n$ is even and $q\in GL_{n}(A)$ is split, with $q^t=q$  then the Chern connection on $GL_n$ attached to $q$ admits a correspondence structure $(\Gamma_p)$ with the following properties: each $\Gamma_p$ has left degree $2^n$; for $n=2$,  $\Gamma_p$  is irreducible, it is partially induced from $E^+$,  and it is not induced from $E^+$; and, for $n=2$,  $(\Gamma_p)$ has $(1,1)$-curvature satisfying $\varphi^*_{p\overline{p}'}\neq0$,  hence  $\varphi_{p\overline{p}'*}\neq 0$, for all $p,p'$. 
\end{theorem}

To prove the theorem we need a preparation. Let us consider an 
$n\times n$ matrix $y$ of indeterminates. Let $B$ be a ring where $2$ is invertible and let ${\mathfrak g}(B)$ be the  $B$-module of all $n\times n$ matrices with entries in $B$, viewed as an associative ring.  Consider the $B$-linear map
$${\mathfrak g}(B)\ra {\mathfrak g}(B),\ \ \ z\mapsto \frac{zb+bz}{2}$$
 (which is the Jordan multiplication by $b$) and let 
$\jor(b)\in B$ be the determinant of this map.

\begin{example}
\label{jordanny}
For $n=2$ we have the following trivially checked formula:
$$\jor(b)=\frac{1}{4}\cdot (\text{tr}(b))^2\cdot \det(b).$$
\end{example}

In what follows, for a matrix $M$ with entries in a ring we denote by $(M)$ the ideal in that ring generated by the entries of $M$. Also, if $B$ is a field, then $dim_B$ below means dimension of a $B$-vector space. Finally let $y$ be a $n\times n$ matrix of indeterminates. In particular we have the polynomial $\jor(y)\in \bZ[1/2][y]$.

Now let $B$ be a ring in which $2$ is invertible, let $b\in GL_n(B)$, and let 
$$C=
B[y,\jor(y)^{-1}]/(y^2-b).$$

\begin{lemma}
\label{lem}

\ 

1) $C$ is \'{e}tale over $B$.

2) If $B$ is a field and $b$ is a scalar matrix then
$\text{dim}_BC=2$.

3) If $B$ is a field and $b$ has distinct eigenvalues then $\text{dim}_BC=2^n$. 

4) If $B$ is a field and  $1$ is not an eigenvalue of $b$ then $\det(y+1)\in C^{\times}$.\end{lemma}

\medskip

{\it Proof of Lemma \ref{lem}}. Assertion 4 is trivial. To check assertion 1  consider a commutative diagram of rings
$$\begin{array}{rcl}
B & \stackrel{\pi}{\longrightarrow} & C\\
v \downarrow & \  & \downarrow u\\
D & \stackrel{\rho}{\longrightarrow} & D/I
\end{array}$$
where $\pi$ is the natural map, $I^2=0$, $v(b)=\beta$, $u(y)=\rho(\gamma)$, $\gamma\in {\mathfrak g}(D)$. So $\gamma^2=\beta+\epsilon$ with $\epsilon$ a matrix with entries in $I$. Since $\rho(\jor(\gamma))$ is invertible in $D/I$ and $I$ is nilpotent it follows that $\jor(\gamma)$ is invertible in $D$. Now $\jor(\gamma)$ is a power of $2$ times the determinant of  the linear system $z\gamma+\gamma z=-\epsilon$ with unknowns the entries of a matrix $z$ hence by Cramer's rule this system has a unique solution $\gamma_1\in {\mathfrak g}(D)$ and this solution is in ${\mathfrak g}(I)$. Note that $(\gamma+\gamma_1)^2=\beta$; so there is a unique homomorphism $U:C\ra D$ such that $U\circ \pi=v$ and $\rho\circ U=u$; it is given by $U(y)=\gamma+\gamma_1$.

To check assertion 2 say $b=\lambda\cdot 1$ with $\lambda\in B^{\times}$. Since, by assertion 1, $C$ is a product of finite field extensions of $B$ it is enough to show that there are exactly $2$ $B$-algebra homomorphisms $C\ra \overline{B}$, where $\overline{B}$ is a fixed algebraic closure of $B$.  Consider such a homomorphism,  $y\mapsto \gamma$, where $\gamma$ is a matrix with coefficients in $\overline{B}$. So $\gamma^2=\lambda \cdot 1$, $\jor(\gamma)\neq 0$. So $\gamma$ is diagonalizable with eigenvalues $\pm \sqrt{\lambda}\in \overline{B}$. It is enough to show that $\gamma$ is a scalar matrix. Let us assume the contrary and derive a contradiction. Indeed, in this case,
$\gamma=u\eta u^{-1}$, $u\in GL_n(\overline{B})$,  
$$\eta=\sqrt{\lambda}\cdot \left( \begin{array}{rrl}
1 & 0 & 0\\
0 & -1 & 0\\
0 & 0 & \eta_{n-2}\end{array}\right),$$
where $\eta_{n-2}$ is a diagonal $(n-2)\times(n-2)$ matrix (with $\pm 1$ on the diagonal).
Since $\jor(\gamma)\neq 0$ the map ${\mathfrak g}(\overline{B})\ra {\mathfrak g}(\overline{B})$, $z\mapsto z\gamma+\gamma z=zu\eta u^{-1}+u\eta u^{-1}z$, is bijective; hence the map
$z\mapsto u^{-1}(z\gamma+\gamma z)u=u^{-1}zu\eta+\eta u^{-1}zu$ is bijective; hence the map
$z\mapsto z\eta+\eta z$ is bijective. But if one takes 
$$z_0=\left( \begin{array}{rrl}
0 & 1 & 0\\
1 & 0 & 0\\
0 & 0 & 0_{n-2}\end{array}\right),$$
one gets $z_0\eta+\eta z_0=0$, a contradiction, which ends the proof of assertion 2.

To check assertion 3 let $\lambda_1,...,\lambda_n\in \overline{B}^{\times}$ be the distinct eigenvalues of $b$, let $\lambda$ be the diagonal matrix with these eigenvalues on the diagonal, and write $b=u\lambda u^{-1}$, with $u\in GL_n(\overline{B})$. Again, it is enough to show that there are exactly $2^n$ $B$-algebra homomorphisms $C\ra \overline{B}$.  Consider such a homomorphism,  $y\mapsto \gamma$, where $\gamma$ is a matrix with coefficients in $\overline{B}$, and $\gamma^2=b^2$. Hence $\gamma$ has distinct eigenvalues, so is diagonalizable, so $\gamma=u\eta u^{-1}$
where $\eta=\text{diag}(\epsilon_1\sqrt{\lambda_1},...,\epsilon_n\sqrt{\lambda_n})$, where $\epsilon_i\in \{\pm 1\}$ and $\sqrt{\lambda_i}$ is a fixed square root to $\lambda_i$. So there are at most $2^n$ homomorphisms $C\ra \overline{B}$. To conclude it is enough to show that for any choice of $\epsilon_1,...,\epsilon_n\in \{\pm 1\}$, and for $\gamma=u\eta u^{-1}$,
$\eta=\text{diag}(\epsilon_1\sqrt{\lambda_1},...,\epsilon_n\sqrt{\lambda_n})$,
 we have $\jor(\gamma)\neq 0$. Assume
$\jor(\gamma)=0$ for some choice of the $\epsilon_i$s and seek a contradiction. 
From $\jor(\gamma)=0$ we know there exists $0\neq z\in {\mathfrak g}(\overline{B})$
such that $\gamma z+z\gamma=0$. Exactly as in the proof of assertion 2 one can find $0\neq w\in {\mathfrak g}(\overline{B})$ such that $\eta w+w\eta=0$, hence 
$\eta w\eta^{-1}=-w$, hence 
$$\frac{\epsilon_i \sqrt{\lambda_i}}{\epsilon_j \sqrt{\lambda_j}}w_{ij}=-w_{ij}$$
for all $i,j$. Choose $i,j$ such that $w_{ij}\neq 0$. Then $i\neq j$ and $\sqrt{\lambda_i}=\pm \sqrt{\lambda_j}$, hence $\lambda_i=\lambda_j$, a contradiction. Assertion 3 is proved.
\qed

\medskip

{\it Proof of Theorem \ref{rational}}.  
Recall that $G=GL_n=Spec\ A[x,\det(x)^{-1}]$ and set $B=A[x,\det(x)^{-1}]$. Let $y$ be a new $n\times n$ matrix of indeterminates.
Define
$$\begin{array}{rcl}
g_p(x) & = & \det(x)\cdot \det(x^{(p)})\cdot \det((x^tqx)^{(p)})  \in  A[x],\\
\  & \   & \  \\
f_p(x) & = & (x^{(p)t}\phi_p(q)x^{(p)})^{-1}(x^tqx)^{(p)}  \in  A[x,g_p(x)^{-1}],\end{array}$$
$$C_p  =  \frac{A[x,g_p(x)^{-1},y, \det(y+1)^{-1}, \jor(y)^{-1}]}{(y^2-f_p(x))},$$
$$Y_{p/G} = Spec\ C_p. $$
Note that for the identity matrix $1$  we have $\jor(1)=1$.
Consider the correspondence $\Gamma_{p/G}=(Y_{p/G},\pi_{p/G},\varphi_{p/G})$ on $G$ where $\pi_{p/G}$ is defined by the natural map (still denoted by) $\pi_{p/G}:B\ra C_p$, $\pi_{p/G}(x)=x$, and $\varphi_{p/G}$ is defined by the map (still denoted by) $\varphi_{p/G}:B\ra C_p$, $\varphi_{p/G}(x)= x^{(p)}y$. Finally define the correspondence 
$\Gamma_p:=\Gamma_{p/G}\otimes E$ on $E$. We will show that $(\Gamma_p)$ possesses  the various required properties in assertions 1, 2, 3. We do this by proving a series of claims as follow.

\medskip

{\it Claim 1}.  $\pi_{p/G}$ is \'{e}tale. 

This follows from Lemma \ref{lem}.
 
 \medskip

{\it Claim 2}. $\pi_{p/G}\otimes {\mathbb F}_p:B\otimes {\mathbb F}_p \ra C_p\otimes {\mathbb F}_p$ is an isomorphism whose inverse sends $y$ into the identity matrix, $1$.

\medskip

Indeed 
$$C_p\otimes {\mathbb F}_p=(A/pA)[x,\det(x)^{-1},y,\det(y+1)^{-1}, \jor(y)^{-1}]/(y^2-1).$$
By Lemma \ref{lem}, ${\mathbb F}_p[y,\jor(y)^{-1}]/(y^2-1)$ is \'{e}tale over ${\mathbb F}_p$ so it is a finite product of finite fields.
In order to prove Claim 2 it is enough to prove that 
$$Z:=Spec\ {\mathbb F}_p[y,\det(y+1)^{-1},
\jor(y)^{-1}]/(y^2-1)
\simeq Spec\  {\mathbb F}_p$$
 under $y\mapsto 1$ so it is enough to prove that the scheme $Z$
  has exactly one $\overline{\mathbb F}_p$-point (which must then be the point given by the identity matrix $y=1$). 
  This follows from Lemma \ref{lem}. Alternatively, 
  let $y\mapsto \gamma\in {\mathfrak g}(\overline{\mathbb F}_p)$ be an $\overline{\mathbb F}_p$-point of $Z$; hence $\gamma^2=1$ and $\det(\gamma+1)\neq 0$. Since the minimal polynomial of $\gamma$ divides $t^2-1$ and $p\neq 2$ the minimal polynomial has distinct roots so $\gamma$ is diagonalizable, 
with eigenvalues $\pm 1$. Since  $-1$ is not an eigenvalue of $\gamma$, we must have  $\gamma=1$.

\medskip

{\it Claim 3}. $\pi_{p/G}^{\widehat{p}}:B^{\widehat{p}}\ra C_p^{\widehat{p}}$
is an isomorphism.

\medskip

Indeed this follows from Claim 2 and from the fact that $p$ is a non-zero divisor in $B$ and $C_p$ (which, in its turn follows from the fact that $C_p$ is \'{e}tale over $B$, cf. Claim 1).

\medskip

{\it Claim 4}. $\varphi_{p/G}^{\widehat{p}}=\pi_{p/G}^{\widehat{p}}\circ \phi_p: B^{\widehat{p}}\ra C_p^{\widehat{p}}$.

\medskip

Indeed, by \ref{turandot}, we have
$\pi_{p/G}^{\widehat{p}}(\phi_p(x))=x^{(p)}f_p(x)^{1/2}$ (where the $1/2$ power is defined by the appropriate $p$-adic series) while $\varphi_{p/G}^{\widehat{p}}(x)=x^{(p)}y$. So
 it is enough to show that $f_p(x)^{1/2}$  and $y$ are equal in $C_p^{\widehat{p}}$. 
Now this follows from the fact that their squares are equal in $C_p^{\widehat{p}}$ and they are both congruent to $1$ in $C_p^{\widehat{p}}$; the fact that $f_p(x)^{1/2}\equiv 1$ mod $p$ is clear from the series expression while the fact that $y\equiv 1$ mod $p$ in $C_p^{\widehat{p}}$ follows from Claim 2.

At this point note that assertion 1  of our Theorem  follows from Claims 1, 3, 4.

Let us check assertion 2 of the Theorem. We first check that if $\Gamma_p=(Y_p,\pi_p,\varphi_p)$ then $Y_p$ is irreducible 
 and $\pi_p$ and $\varphi_p$ have degrees $2$ and $2p^4$ respectively. If $$x=\left(\begin{array}{cc}a & b\\ c & d\end{array}\right),\ \ \  \beta_p=\frac{\det(x)^p}{\det(x^{(p)})}=\frac{(ad-bc)^p}{a^pd^p-b^pc^p}$$
 then
$$Y_p=Y_{p/G}\otimes E=Spec\ \frac{E[y,\det(y+1)^{-1},\jor(y)^{-1}]}{(y^2-\beta_p\cdot 1)}.$$
Let $t_p=t$ be an indeterminate and set 
$$F_p=\frac{E[t]}{(t^2-\beta_p)}.$$ Since
$\beta_p^{-1}$ is the inverse of a square  times a product of coprime linear forms in ${\mathbb Q}[a,b,c,d]$ it follows that $t^2-\beta_p$ is irreducible in $E[t]$, hence $F_p$ is a field of degree $2$ over $E$. We claim that
$Y_p\simeq Spec\ F_p$.
Indeed the map 
$$E[y]\ra F_p,\ \ y\mapsto \left(\begin{array}{cc}t & 0\\ 0 & t\end{array}\right)$$
sends $(y^2-\beta_p\cdot 1)$ into $0$ and induces a map 
$$\cO(Y_p)\ra F_p$$
because $\jor \left( \begin{array}{cc}t & 0\\ 0 & t\end{array}\right)=t^{4}\in F_p^{\times}$
and $\det \left(\begin{array}{cc}t+1 & 0\\ 0 & t+1\end{array}\right)\in F_p^{\times}$. Since $Y_p$ is the disjoint union of spectra of  fields  it follows that  $Spec\ F_p$ is a connected component of $Y_p$. In order to show there are no other connected components it is enough to show that $\cO(Y_p)$ has dimension at most  two 
over $E$; this follows from Lemma \ref{lem}. This ends the proof of the fact that $Y_p$ is irreducible and of the fact that $\pi_p$ has degree $2$. To prove that $\varphi_p$ has degree $2p^4$ we need to show that, for $\theta:=\theta_p\in F_p$ the class of $t=t_p$, the degree of the extension
\begin{equation}
\label{extension}
K(a^p\theta,b^p\theta,c^p\theta,d^p\theta)\subset K(a,b,c,d,\theta)\end{equation}
is $2p^4$. Set $u=a/d$, $v=b/d$, $w=c/d$, $M=K(u,v,w)$, and note that $\theta^2=\beta_p\in M$. Then the extension \ref{extension} is the composition of the following extensions:
$$
\begin{array}{rcll}
K(a^p\theta,b^p\theta,c^p\theta,d^p\theta) & = & 
K(u^p,v^p,w^p,d^p\theta) & \ \\
\  & \subset & M(d^p\theta) & (\star)\\
\  & \subset & M(\theta,d^p) & (\star\star)\\
\  & \subset & M(\theta,d) & (\star\star\star)\\
\  & = & K(a,b,c,d,\theta). &\ 
\end{array}
$$
Now the extensions $(\star)$ and $(\star\star\star)$ obviously have degrees $p^3$ and $p$ respectively. We claim that the extension $(\star\star)$ has degree $2$. It has degree $1$ or $2$ so we need to show that
$\theta\not\in M(d^p\theta)$. Assume the contrary. Then, writing $z$ for the indeterminate $d^p$ over $M$ we have $\theta\cdot F(\theta z)=G(\theta z)$ for two polynomials $F,G$ with coefficients in $M$.
Set $z=\theta^{2i-1}$ for $i=1,2,3,...$. We get $\theta\cdot F(\beta_p^i)=G(\beta_p^i)$ for all $i$. Since $\beta_p$ is not a root of unity we have $F(\beta_p^i)\neq 0$ for some $i$; hence we get $\theta\in M$, a contradiction. This ends the proof that $\varphi_p$ has degree $2p^4$. 

Next we show that $\varphi_{pp'}^*=0$ i.e. we show that the $(\Gamma^*_p)$'s commute. View $Y_p=Spec\ F_p$ and view $\pi_p:E\ra F_p=E+\theta_p E$ as the natural inclusion, $\theta_p^2=\beta_p$. So $\text{tr}_{\pi_p}(\theta_p)=0$. Also $\varphi_p(x)=\theta_p x^{(p)}$, i.e. $\varphi_p(a)=\theta_pa^p$, $\varphi_p(b)=\theta_pb^p$, etc.
Recall the spaces $E^+$ and $E^-$ in \ref{E+}. Note that $\beta_p\in E^+$ and $\varphi_p(E^+)\subset E^+$, hence 
 \begin{equation}
 \label{regularity}
 \Gamma^*_p(E^+)\subset E^+,\ \ \Gamma^*_{|E^+}=2\varphi_p.
 \end{equation}

 {\it Claim 5}.  $\varphi_p\circ \varphi_{p'}$ and $\varphi_{p'}\circ \varphi_p$ have the same value on any monomial of degree $2$ in  $\{a,b,c,d\}$.

Indeed
 let $a^ib^jc^kd^l$ be such a monomial, $i+j+k+l=2$.  
Then
$$\begin{array}{rcl}
 \Gamma^*_{p'}\Gamma^*_p(a^ib^jc^kd^l) & = &
  \Gamma^*_{p'}\left(\text{tr}_{\pi_p}\left(a^{pi}b^{pj}c^{pk}d^{pl}\frac{\det(x)^p}{\det(x^{(p)})}\right)\right)\\
 \  & \ & \  \\
 \  & = & 2  \Gamma^*_{p'}\left(
 a^{pi}b^{pj}c^{pk}d^{pl}
 \frac{\det(x)^p}{\det(x^{(p)})}\right)\\
 \  & \ & \  \\
 \  & \  & \  \\
 \  & = & 4 a^{pp'i}b^{pp'j}c^{pp'k}d^{pp'l}\cdot
 \left(\frac{\det(x)^{p'}}{\det(x^{(p')})}\right)^p\cdot
 \frac{\det(x^{(p')})^p}{\det(x^{(pp')})}\\
 \  & \  & \  \\
 \  & = & 4 a^{pp'i}b^{pp'j}c^{pp'k}d^{pp'l}\cdot
 \frac{\det(x)^{pp'}}{\det(x^{(pp')})}\\
 \  & \  & \  \\
 \  & = &  \Gamma^*_{p}\Gamma^*_{p'}(a^ib^jc^kd^l). 
 \end{array}
  $$
  From Claim 5 we get that 
  $\varphi_p\circ \varphi_{p'}$ and $\varphi_{p'}\circ \varphi_p$
   are equal as maps $E^+\ra E^+$. Hence $\Gamma^*_p$ and $\Gamma^*_{p'}$ commute on $E^+$.
On the other hand 
\begin{equation}
\label{regularity1}
\Gamma^*_p(E^-)\in \text{tr}_{\pi_p}(\varphi_p(a)\varphi_p(E^+)) \subset \text{tr}_{\pi_p}(\theta_p E)=0.\end{equation}
 Similarly $\Gamma^*_{p'}$ vanishes on $E^-$. So $\Gamma_p^*$ and $\Gamma^*_{p'}$ commute on $E^-$ as well, so $\Gamma_p^*$ and $\Gamma^*_{p'}$ commute on $E$.
 
 Note, at this point, that \ref{regularity} and \ref{regularity1} show that $\Gamma_p$ is induced from $E^+$. 
 
 We next show that $\varphi_{p\overline{p}'}\neq 0$. Using the fact that $\Gamma^*_{\overline{p}'}(x)=x^{(p')}$, we have:
 $$
 \begin{array}{rcl}
 \Gamma^*_{\overline{p}'}\Gamma^*_p(a^2) & = &
 \Gamma^*_{\overline{p}'}\left(\text{tr}_{\pi_p}\left(a^{2p}\frac{\det(x)^p}{\det(x^{(p)})}\right)\right)\\
 \  & \ & \  \\
 \  & = & 2 a^{2pp'} \frac{\det(x^{(p')})^p}{\det(x^{(pp')})},\\
 \  & \  & \  \\
 \Gamma_p^*\Gamma_{\overline{p}'}^* (a^2) & = & 
 \Gamma^*_p(a^{2p'})\\
 \  & \  & \  \\
 \  & = & 2a^{2pp'} \left(\frac{\det(x)^p}{\det(x^{(p)})}\right)^{p'}.
 \end{array}
 $$
So $\Gamma^*_{\overline{p}'}\Gamma^*_p\neq \Gamma_p^*\Gamma_{\overline{p}'}^*$ 
because $\frac{\det(x^{(p')})^p}{\det(x^{(pp')})}\neq \left(\frac{\det(x)^p}{\det(x^{(p)})}\right)^{p'}$; the latter can be checked by taking $a=c=d=1$, clearing out the denominators and picking out the coefficient of $b$. 

In order to end the proof of assertion 2 we need to show that $\Gamma_p$ and $\Gamma_{p'}$ partially commute, but do not commute. Recall that we established an identification $Y_p=Y_{\Gamma_p}=Spec\ F_p$, $Y_{p'}=Y_{\Gamma_{p'}}=Spec\ F_{p'}$ where 
$$F_p=\frac{E[t_p]}{(t^2_p-\beta_p)}, \ \ F_{p'}=\frac{E[t_{p'}]}{(t^2_{p'}-\beta_{p'})},$$
where $\pi_p,\pi_{p'}$ identify with the inclusions of $E$ into $F_p$ and $F_{p'}$ respectively and $\varphi_p:E\ra F_p$, $\varphi_{p'}:E\ra F_{p'}$ are given by
$\varphi_p(x)=t_px^{(p)}$ and $\varphi_{p'}(x)=t_{p'}x^{(p')}$. Then we have $Y_{\Gamma_{p}\circ \Gamma_{p'}}=Spec\ F_{pp'}$, $Y_{\Gamma_{p'}\circ \Gamma_p}=Spec\ F_{p'p}$,  where
$$F_{pp'}=\frac{E[t_p,t_{p'}]}{(t^2_{p'}-\beta_{p'},t_{p}^2-\varphi_{p'}(\beta_{p}))},\ \ F_{p'p}=\frac{E[t_p,t_{p'}]}{(t^2_p-\beta_p,t_{p'}^2-\varphi_p(\beta_{p'}))}.$$
We need to prove that the two maps
\begin{equation}
\label{twomaps}
Spec(F_{pp'}\otimes_E F_{p'p})\ra Spec\ E\end{equation}
defined by the ring maps
$$E\ra F_{pp'}\otimes_E F_{p'p},\ \  x\mapsto 1\otimes  t_{p'}t_p^{p'}x^{(pp')}\ \ \ \text{and}\ \ \ x\mapsto  t_p t_{p'}^px^{(pp')}\otimes 1,$$
do not coincide on $Spec(F_{pp'}\otimes_E F_{p'p})$ but {\it do} coincide on a connected component of the latter. Note that $x^{(pp')}$ has entries in $E$ so may be balanced left and right of the $\otimes$ sign. Set 
$$\psi^{\pm}:=t_{p}t_{p'}^{p}\otimes 1 \pm 1\otimes t_{p'} t_{p}^{p'} \in F_{pp'}\otimes_E F_{p'p}.$$
We have that $\psi^+\neq 0$ and $\psi^-\neq 0$ (otherwise we get a contradiction with the fact that $F_{pp'}\otimes_E F_{p'p}$ has dimension $16$ over $E$). On the other hand we have
$$\psi^+\psi^-=t^2_{p}t_{p'}^{2p}\otimes 1 -1\otimes t^2_{p'} t_{p}^{2p'} 
=\varphi_{p'}(\beta_p)\beta_{p'}^p-\varphi_{p}(\beta_{p'})\beta_{p}^{p'}.
$$
Now an easy computation gives
$$\varphi_{p'}(\beta_p)\beta_{p'}^p=\frac{(\det(x))^{pp'}}{\det(x^{(pp')})}=\varphi_{p}(\beta_{p'})\beta_{p}^{p'},$$
hence $\psi^+ \psi^-=0$. So neither $\psi^+$ nor $\psi^-$ is invertible in $F_{pp'}\otimes_E F_{p'p}$. Let $M^+$ and $M^-$ be maximal ideals in $F_{pp'}\otimes_E F_{p'p}$ containing $\psi^+$ and  $\psi^-$ respectively. Then the two maps in \ref{twomaps} coincide on $Spec((F_{pp'}\otimes_E F_{p'p})/M^-)$ and {\it do not} coincide on $Spec((F_{pp'}\otimes_E F_{p'p})/M^+)$. This ends the proof of assertion 2.

Let us prove  assertion 3 of the theorem.

We start by proving that for $n$ even and $q$ split with $q^t=q$, $\Gamma_p$ has left degree $2^n$. To check this it is enough to check that $f_p(x)$ has distinct eigenvalues in an algebraic closure $\overline{E}$ of $E$ and none of these eigenvalues is $1$; cf. Lemma \ref{lem}. Let us introduce some notation. For any square matrix $M$ with entries in a ring we denote by  $Char(M)$ the characteristic polynomial
of $M$; furthermore, for any monic polynomial $f$ with coefficients in our  ring we denote by $Dis(f)$ the discriminant of $f$; so we may consider the element
$Dis(Char(M))$ in our ring. In particular we may consider the element
  ${\mathcal D}_p=Dis(Char(f_p(x)))\in A[x,g_p(x)^{-1}]$. 
Let $n=2r$ and write $q=q_n$. Then
$q_n=w^{-1}\tilde{q}_nw$ where $w$ is a permutation matrix and 
$$\tilde{q}_n=diag(q_2,...,q_2),\ \ \ \ \ q_2=\left(\begin{array}{cc} 0 & 1 \\ 1 & 0\end{array}\right)$$
Since $w$ is a permutation matrix we have $w^{-1}=w^t$ and $(wx)^{(p)}=wx^{(p)}$ so 
$$x^{(p)t}q_nx^{(p)}=x^{(p)t}w^t\tilde{q}_nwx^{(p)}=(wx)^{(p)t}\tilde{q}_n(wx)^{(p)},$$
$$x^{t}q_nx=x^{t}w^t\tilde{q}_nwx=(wx)^{t}\tilde{q}_n(wx),$$
and hence, for $\tilde{x}=wx$, we have
\begin{equation}
\label{vint}
f_p(x)=(\tilde{x}^{(p)t}\tilde{q}_n\tilde{x}^{(p)})^{-1}(\tilde{x}^t\tilde{q}_n\tilde{x})^{(p)}.
\end{equation}
For $1\leq l,m\leq r$ consider the matrices
$$z_{lm}=\left(\begin{array}{ll}
\tilde{x}_{2l-1,2m-1} & \tilde{x}_{2l-1,2m} \\ \tilde{x}_{2l,2m-1} & \tilde{x}_{2l,2m} \end{array}\right).$$
Let $\cC_p:=\cC(\tilde{x})$ be  the characteristic polynomial of   \ref{vint}.
Then if one sets $z_{lm}=0$ in $\cC_p(\tilde{x})$ for all $l\neq m $ one obtains the product
$\cC_p^1(z_{11})...\cC_p^r(z_{rr})$ where $\cC_p^l(z_{ll})$ is the characteristic polynomial of
$$f_p^l:=(z_{ll}^{(p)t}q_2z_{ll}^{(p)})^{-1}(z_{ll}^tq_2z_{ll})^{(p)}=:\left(\begin{array}{cc} u_l & v_l\\ w_l& u_l\end{array}\right).$$
Set
$$z_{ll}=\left(\begin{array}{cc} a & b\\ c& d\end{array}\right).$$
Then, by a direct computation, using \ref{turandot},  we have $u_l=u$, $v_l=v$, $w_l=w$, where
\begin{equation}
\label{uvw}\begin{array}{rcl}
 u & = & \frac{(a^pd^p+b^pc^p)(ad+bc)^p-2^{p+1}a^pb^pc^pd^p}{(a^pd^p-b^pc^p)^2},\\
 \  & \  & \  \\
 v & = & b^pd^p\cdot \frac{2^p(a^pd^p+b^pc^p)-2(ad+bc)^p}{(a^pd^p-b^pc^p)^2},\\
 \  & \  & \  \\
 w & = & a^pc^p\cdot \frac{2^p(a^pd^p+b^pc^p)-2(ad+bc)^p}{(a^pd^p-b^pc^p)^2}.\end{array}
 \end{equation}
  Note that $v_l,w_l\neq 0$; hence the discriminant of the characteristic polynomial of $f_p^l$, which equals  $4v_lw_l$, is non-zero.
One then easily gets that the discriminant of $\cC_p^1(z_{11})...\cC_p^r(z_{rr})$ is non-zero 
in 
$$A[z,g_p(z)^{-1}]:=A[z_{11},...,z_{rr},g_p(z_{11})^{-1},...,g_p(z_{rr})^{-1}],$$
where the latter $g_p(z_{ii})$ are the corresponding quantities attached to $q_2$. Let $\tilde{g}_p$ be the  polynomial analogous to $g_p$ constructed for $\tilde{q}_n$ rather than $q_n$.
 We conclude that ${\mathcal D}_p\neq 0$  in $A[x,g_p(x)^{-1}]$ because the image of ${\mathcal D}_p$ under the map 
$$A[x,g_p(x)^{-1}]=A[\tilde{x},\tilde{g}_p(\tilde{x})^{-1}]\ra A[z,g_p(z)^{-1}]$$
$$z_{lm}\mapsto 0\ \ \text{for}\ \ l\neq m,$$
$$z_{ll}\mapsto z_{ll}$$
is $\neq 0$. Finally to show that $1$ is not an eigenvalue for $f_p(x)$  it is enough to show that $1$ is not an eigenvalue for any of the matrices $f_p^l$. Assume $1$ is an eigenvalue of $f_p^l$ and seek a contradiction. Our assumption implies that  $(u_l-1)^2-v_lw_l=0$, hence $v_lw_l$ is a square in $K(a,b,c,d)$, where $K$ is the fraction field of $A$. But $v_lw_l$ is a square in 
$K(a,b,c,d)$
times $a^pb^pc^pd^p$, hence $a^pb^pc^pd^p$ is a square in $K(a,b,c,d)$, a contradiction.

Next we assume $n=2$, hence 
$q=\left(\begin{array}{cc} 0 & 1 \\ 1 & 0\end{array}\right)$,
 and we prove the remaining claims in assertion 3 of the Theorem.
Write $x=\left(\begin{array}{cc}a & b\\ c & d\end{array}\right)$.
Then recall from the argument above that
$f_p=\left(\begin{array}{cc}u & v\\ w & u\end{array}\right)$ 
where $u,v,w$ are given by \ref{uvw}. Recall also that, by that argument, $1$ is not an eigenvalue of $f_p$ and $vw\neq 0$. So by Lemma \ref{lem} and the proof of assertion 1 of the theorem we have $Y_p=Spec\ F_p$ where
$$F_p=\frac{E[y,\jor(y)^{-1}]}{(y^2-f_p)}.$$
If $\left(\begin{array}{cc}\alpha & \beta\\ \gamma & \d\end{array}\right)$ is the image of $y$ in the $2\times 2$ matrices with entries in $F_p$ then
$\pi_p:E\ra F_p$ identifies with the natural inclusion, and 
\begin{equation}
\label{eqq}
\begin{array}{ll}
\varphi_p(a)  = a^p\alpha+b^p\gamma, & \varphi_p(b) =  a^p\beta+b^p\d,\\
\  & \ \\
\varphi_p(c)  =  c^p\alpha+d^p\gamma, & \varphi_p(d)  =  c^p\beta+d^p\d.\end{array}\end{equation}
We claim that, if $t$ is one variable then we have  natural isomorphisms
\begin{equation}
\label{damn}
\frac{E[t]}{(t^4-4ut^2+4vw)}\simeq \frac{E[y]}{(y^2-f_p)}\simeq F_p,\ \ \ t\mapsto \tau:=\alpha+\d.\end{equation}
To prove \ref{damn}  we will first prove the following:

{\it Claim 6}.  The map $h:E[t]\ra  \frac{E[y]}{(y^2-f_p)}$, $t\mapsto \tau$, sends $t^4-4ut^2+4vw$ into $0$.

{\it Claim 7}. $\tau$ is invertible in $\frac{E[y]}{(y^2-f_p)}$.

{\it Claim 8}.  The map $h$ is surjective. 

\noindent To check Claim 6, note that 
\begin{equation}
\label{tool}
\alpha^2+\beta\gamma=\beta\gamma+\d^2=u,\ \ (\alpha+\d)\beta=v,\ \ \ (\alpha+\d)\gamma=w.\end{equation}
Set  $\eta=\alpha-\d$. Then 
\begin{equation}
\label{force}
\tau\eta=\alpha^2-\d^2=0,\ \ \ \tau^2+\eta^2=2(\alpha^2+\d^2).\end{equation}
By \ref{tool} $$\alpha^2+\d^2+2\beta\gamma=2u,$$ hence, by \ref{force},  
$$\tau^2+\eta^2+4\beta\gamma=4u,$$ hence $$\tau^4+\tau^2\eta^2+4\beta\gamma\tau^2=4u\tau^2,$$ hence by \ref{force} and \ref{tool}, 
\begin{equation}
\label{unamica}
\tau^4-4u\tau^2+4vw=0.
\end{equation}
 This ends the proof of Claim 6.
Note that by Claim 6 we have an induced homomorphism
$$\overline{h}:\frac{E[t]}{(t^4-4ut^2+4vw)}\ra \frac{E[y]}{(y^2-f_p)}.$$

Claim 7 follows from the fact that the class of $t$ in $\frac{E[t]}{(t^4-4ut^2+4vw)}$ is invertible (because $vw\neq 0$) and $\tau$ is the image of this class via the homomorphism $\overline{h}$.

Claim 8 is checked as follows. By \ref{tool} and Claim 7, 
\begin{equation}
\label{z1}
\beta=\tau^{-1}v,\ \ \gamma=\tau^{-1}w
\end{equation}
 belong to the image of the map $h$. Also, since by \ref{force}, $\tau\eta=0$, we get $\eta=0$ hence 
 \begin{equation}
 \label{z2}
 \alpha=\d=\tau/2,
 \end{equation}
 so $\alpha$ and $\beta$ belong to the image of the map $h$. 

Now \ref{damn} can be checked as follows. By Claim 7 and by Example \ref{jordanny} the image of $\jor(y)$ in $\frac{E[y]}{(y^2-f_p)}$ is invertible. So
 $F_p\simeq\frac{E[y]}{(y^2-f_p)}$.
 Since, by Theorem \ref{rational}, $F_p$ has dimension $4$ over $E$ and $\overline{h}$ is surjective (cf. Claim 8), \ref{damn} follows.
 
 Note that by \ref{damn} we have 
 \begin{equation}
 \label{Trr}
 \text{tr}_{\pi_p}(\tau)=\text{tr}_{\pi_p}(\tau^{-1})=0.
 \end{equation}
 
 Recall the automorphism $\iota: E\ra E=K(x)$, $\iota(x)=- x$ and the spaces
$E^+,E^-$ in \ref{E+}. Note that 
\begin{equation}
\label{inE+}
u,v,w\in E^+,\ \ \end{equation}
hence
\begin{equation}
\label{inE++}
\text{tr}_{\pi_p}(\tau^2)=8u\in E^+,\ \  \text{tr}_{\pi_p}(\tau^{-2})=\frac{2u}{vw}\in E^+.\end{equation}
Using \ref{eqq}, \ref{z1}, \ref{z2},  \ref{inE+},  one 
gets that $\varphi_p(E^+)\subset E^+ + E^+\tau^2$ and $\varphi_p(E^-)\subset E^+\tau+E^+\tau^{-1}$.
Hence by
\ref{Trr}, \ref{inE++} one gets that
the map $\Gamma_p^*=\text{tr}_{\pi_p}\circ \varphi_p:E\ra F_p\ra E$ satisfies
\begin{equation}
\label{regularity2}
\Gamma^*_p(E^+)\subset E^+,\ \ \Gamma_p^*(E^-)=0,\end{equation}
and  partial induction from $E^+$ in assertion 3 is proved.

In what follows we show that $\varphi^*_{p\overline{p}'}\neq 0$ and that $\Gamma^*_p$ is not induced from $E^+$. Note  that 
\begin{equation}
\label{burta}
\begin{array}{rcll}
\varphi_p(ab) & = & a^{2p}\alpha\beta +a^pb^p(\alpha\d+\beta\gamma)+b^{2p}\gamma\d & \text{(by \ref{eqq}})\\
\  & \  & \ & \\
\  & = & a^{2p}\frac{v}{2}+b^{2p}\frac{w}{2}+a^pb^p(\frac{\tau^2}{4}+vw\tau^{-2}) & \text{(by \ref{z1}, \ref{z2})}\\
\  & \  & \ & \\
\  & = & a^{2p}\frac{v}{2}+b^{2p}\frac{w}{2}+ua^pb^p & \text{(by \ref{unamica})}\\
\  & \  & \ & \\
\  & = & 2^{p-1}a^pb^p & \text{(by \ref{uvw})} \end{array}
\end{equation}
One gets
$$\begin{array}{rcll}
\Gamma^*_{\overline{p}'}\Gamma^*_p(ab) & = & \Gamma_{\overline{p}'}^*(\text{tr}_{\pi_p}(\varphi_p(ab))) & \\
\  & \  & \ & \\
\  & = & \Gamma_{\overline{p}'}^*(\text{tr}_{\pi_p}(2^{p-1}a^pb^p)) & \text{(by \ref{burta})}\\
\  & \  & \ & \\
\  & =  & \Gamma_{\overline{p}'}^*(2^{p+1}a^pb^p) & \text{(since $deg(\pi_p)=4$)}\\
\  & \  & \ & \\
\  & = & 2^{p+1}a^{pp'}b^{pp'}. &\end{array}
$$
On the other hand, similarly,  we have
$$
\begin{array}{rcll}
\Gamma^*_p \Gamma^*_{\overline{p}'}(ab) & = & \Gamma^*_p(a^{p'}b^{p'}) &\\
\  & \  &  \  & \\
\  & = & \text{tr}_{\pi_p}((\varphi_p(ab))^{p'}) & \\
\  & \  & \ & \\
\  & = & \text{tr}_{\pi_p}((2^{p-1}a^pb^p)^{p'}) & \\
\  & \  & \ & \\
\  &  = & 4\cdot 2^{(p-1)p'}a^{pp'}b^{pp'}.\end{array}
$$
So we get 
$$
\varphi^*_{p\overline{p}'}(ab)=\frac{2^{p+1}(1-2^{(p-1)(p'-1)})}{pp'}a^{p^2}b^{p^2}\neq 0,\ \ \text{for $p\neq p'$, and}$$
$$
\varphi^*_{p\overline{p}}(ab)=\frac{2^{p+1}(1-2^{(p-1)^2})}{p}a^{p^2}b^{p^2}\neq 0,$$
so $\varphi^*_{p\overline{p}'}\ne0$ for all $p,p'$. 
Finally a direct computation gives
$$\varphi_p(ad)=(a^pc^p\frac{v}{2}+b^pd^p\frac{w}{2}+ub^pc^p)+(a^pd^p-b^pc^p)\frac{\tau^2}{4}\not\in E^+,$$
hence $\Gamma^*_p$ is not induced from $E^+$.

To close the proof of assertion 3 we need to show that $F_p$ is a field. By \ref{damn} it is enough to show the following:

{\it Claim 9}. The $E$-subalgebra $E[\tau^2]$ of $F_p$ generated by $\tau^2$ is a field.

{\it Claim 10.} $\tau^2$ is not a square in the field $E[\tau^2]$.

We proceed to proving Claim 9. Note  that the discriminant of the polynomial $y^2-4uy+4vw$,  satisfied by $\tau^2$, has the following simple form:
\begin{equation}
\label{jerry}
16u^2-16vw=16\cdot\frac{(ad+bc)^{2p}-4^pa^pb^pc^pd^p}{(a^pd^p-b^pc^p)^2}.
\end{equation}
We claim that this discriminant is not a square in the field $E=K(a,b,c,d)$. We need to show that $Q:=(ad+bc)^{2p}-4^pa^pb^pc^pd^p$ is not a square in $K(a,b,c,d)$. Assume the contrary. Then $Q=P^2$  with $P\in K[a,b,c,d]$ because $K[a,b,c,d]$ is factorial.
Consider the ${\mathbb G}_m\times {\mathbb G}_m$-action $\star$ on $K[a,b,c,d]$,
${\mathbb G}_m=Spec\ K[t,t^{-1}]$,  with points $g:=(\lambda,\mu)$ of 
${\mathbb G}_m\times {\mathbb G}_m$
acting via $g\star a=\lambda a$, $g\star b=\mu b$, $g\star c=\mu^{-1} c$, $g\star d=\lambda^{-1}d$, respectively. The ring of invariants of this action is $K[ad,bc]$. Since $Q$ is invariant we have $g\star P=\epsilon(g)\cdot  P$, $\epsilon(g)\in \{\pm1\}$,  for all $g$. Since 
${\mathbb G}_m\times {\mathbb G}_m$
 is connected we must have $\epsilon(g)=1$ for all $g$, hence $P(a,b,c,d)=G(ad,bc)$
for some polynomial $G$ in $2$ variables with $K$-coefficients.
Hence, for $s$ a variable, we have 
$$f(s):=Q(s,1,1,1)= (s+1)^{2p}-4^ps^p=G(s,1)^2$$ in $K[s]$, so all the complex roots of $f$ are multiple. Now,
if $f'(s)=df/ds$, we get
$$(s+1)f'(s)-2pf(s)=4^p ps^{p-1}(s-1).$$
So the only root of $f$ is $1$, hence $f(s)=(s-1)^{2p}$, a contradiction, and our Claim  9 is proved.

Next we proceed to proving Claim 10. We have
$\tau^2=2u\pm \sqrt{u^2-vw}$. Assume $\tau^2$ is a square in $E[\tau^2]$. Using \ref{jerry}, and setting $D:=(ad+bc)^{2p}-4^pa^pb^pc^pd^p$, 
we get that
$$2(a^pd^p+b^pc^p)(ad+bc)^p-2^{p+2}a^pb^pc^pd^p\pm 2(a^pd^p-b^pc^p)\sqrt{D}=
(X+Y\sqrt{D})^2$$
for some $X,Y\in E$. Set $A=ad$, $B=bc$. We get
$$
\begin{array}{rcl}
2(A^p+B^p)(A+B)^p-2^{p+2}A^pB^p & = & X^2+DY^2\\
\  & \  & \  \\
2(A^p-B^p) & = & 2XY,
\end{array}
$$
for some $X,Y\in E$. Again , consider the ${\mathbb G}_m\times {\mathbb G}_m$-action $\star$ on $E=K(a,b,c,d)$,
${\mathbb G}_m=Spec\ K[t,t^{-1}]$,  with points $g:=(\lambda,\mu)$ of 
the torus ${\mathbb G}_m\times {\mathbb G}_m$
acting via $g\star a=\lambda a$, $g\star b=\mu b$, $g\star c=\mu^{-1} c$, $g\star d=\lambda^{-1}d$, respectively. The invariant field of this action is $K(A,B)$. By connectivity of the torus we get that $X$ and $Y$ are invariant and hence belong to $K(A,B)$. Set $s=B/A$ and $Z=X/A$. We get
$$
\begin{array}{rcl}
M:=2(1+s^p)(1+s)^p-2^{p+2}s^p & = & Z^2+((1+s)^{2p}-4^ps^p)Y^2=Z^2+PY^2\\
\  & \  & \  \\
N:=2(1-s^p) & = & 2ZY.
\end{array}
$$
Consider the action of ${\mathbb G}_m$ on $K(A,B)$ with $\lambda$ acting on $A,B$ by sending them into $\lambda A, \lambda B$ respectively. Since ${\mathbb G}_m$ is connected and $s$ is invariant it follows that $Z,Y$ are invariant hence belong to $K(s)$. Eliminating $Z$ from the above system we get
$PY^4-MY^2+N^2/4=0$. Since $Y\in K(s)$ we get that the discriminant $M^2-PN^2$ of the latter equation is a square in $K(s)$. In particular $M^2-PN^2$ has even $s$-adic valuation.  On the other hand we have
$$
\begin{array}{rcl}
4^{-1}(M^2-PN^2) & = & ((1+s^p)(1+s)^p-2^{p+1}s^p)^2 -(1-s^p)^2((1+s)^{2p}-4^ps^p)\\
\  & \  & \  \\
\  & \equiv  & (4^p-2^{p+2}+4)s^p\ \ \text{mod $s^{p+1}$}.
\end{array}
$$
Since 
$4^p-2^{p+2}+4\neq 0$ the $s$-adic valuation of $M^2-PN^2$ is $p$ which is odd, a contradiction. This ends the proof of Claim 10, hence the proof of assertion 3, and hence the proof of the Theorem.
\qed

\end{document}